\input amstex
\input xy
\xyoption{all}
\documentstyle{amsppt}
\magnification=\magstep 1
\document
\vsize6.7in

%%%%%%%%%%%%%%%%%%%%%%%%%%%%%%%%%%%%%%%%%%%%%%%%%%%%%
% One can jump to the label : THE PAPER STARTS HERE %
%%%%%%%%%%%%%%%%%%%%%%%%%%%%%%%%%%%%%%%%%%%%%%%%%%%%%

\chardef\oldatsign=\catcode`\@
\catcode`\@=11
\newif\ifdraftmode			% New- print marginal notes if true
\global\draftmodefalse

% End of \draft.


%===============================================================================
% Font loading. Normal text size is 12pt.
%
\font@\twelverm=cmr12 % roman text 
\font@\twelvei=cmmi12 \skewchar\twelvei='177 % math italic
\font@\twelvesy=cmsy10 scaled\magstep1 \skewchar\twelvesy='060 % math symbols
\font@\twelveex=cmex10 scaled\magstep1 % math extensions
\font@\twelvemsa=msam10 scaled\magstep1 % AMS extra symbols
\font@\twelvemsb=msbm10 scaled\magstep1 % AMS extra symbols
\font@\twelvebf=cmbx12 % boldface extended
\font@\twelvett=cmtt12 % typewriter
\font@\twelvesl=cmsl12 % slanted roman
\font@\twelveit=cmti12 % text italic
\font@\twelvesmc=cmcsc10 scaled\magstep1 % caps and small caps
%
%     the normal script font is nine point
%
\font@\ninerm=cmr9 % roman text
\font@\ninei=cmmi9 \skewchar\ninei='177 % math italic
\font@\ninesy=cmsy9 \skewchar\ninesy='60 % math symbols
\font@\ninemsa=msam9
\font@\ninemsb=msbm9
\font@\ninebf=cmbx9
%
%	the normal scriptscript font is 7 point; already loaded by amstex
%	or plain.
%
%     some specialty fonts
%
\font@\ttlrm=cmbx12 scaled \magstep2 % a title font
\font@\ttlsy=cmsy10 scaled \magstep3 % for \AmSTeX in titles
\font@\tensmc=cmcsc10 % caps and small caps
 % a little caps type
 % a little tt type
%
%     twelve point is used for normal printing.
%
\def\normaltype{%   twelve point stuff
	\def\pointsize@{12}%
	\abovedisplayskip18\p@ plus5\p@ minus9\p@
	\belowdisplayskip18\p@ plus5\p@ minus9\p@
	\abovedisplayshortskip1\p@ plus3\p@
	\belowdisplayshortskip9\p@ plus3\p@ minus4\p@
	\textonlyfont@\rm\twelverm
	\textonlyfont@\it\twelveit
	\textonlyfont@\sl\twelvesl
	\textonlyfont@\bf\twelvebf
	\textonlyfont@\smc\twelvesmc
	\ifsyntax@
		\def\big##1{{\hbox{$\left##1\right.$}}}%
	\else
		\let\big\twelvebig@
 \textfont0=\twelverm \scriptfont0=\ninerm \scriptscriptfont0=\sevenrm
 \textfont1=\twelvei  \scriptfont1=\ninei  \scriptscriptfont1=\seveni
 \textfont2=\twelvesy \scriptfont2=\ninesy \scriptscriptfont2=\sevensy
 \textfont3=\twelveex \scriptfont3=\twelveex  \scriptscriptfont3=\twelveex
 \textfont\itfam=\twelveit \def\it{\fam\itfam\twelveit}%
 \textfont\slfam=\twelvesl \def\sl{\fam\slfam\twelvesl}%
 \textfont\bffam=\twelvebf \def\bf{\fam\bffam\twelvebf}%
 \scriptfont\bffam=\ninebf \scriptscriptfont\bffam=\sevenbf
 \textfont\ttfam=\twelvett \def\tt{\fam\ttfam\twelvett}%
 \textfont\msafam=\twelvemsa \scriptfont\msafam=\ninemsa
 \scriptscriptfont\msafam=\sevenmsa
 \textfont\msbfam=\twelvemsb \scriptfont\msbfam=\ninemsb
 \scriptscriptfont\msbfam=\sevenmsb
	\fi
 \normalbaselineskip=\twelvebaselineskip
 \setbox\strutbox=\hbox{\vrule height12\p@ depth6\p@
      width0\p@}%
 \normalbaselines\rm \ex@=.2326ex%
}% End of \normaltype.
%
%
%     ten point is used for smalltype
%
\def\smalltype{%   ten point stuff
	\def\pointsize@{10}%
	\abovedisplayskip12\p@ plus3\p@ minus9\p@
	\belowdisplayskip12\p@ plus3\p@ minus9\p@
	\abovedisplayshortskip\z@ plus3\p@
	\belowdisplayshortskip7\p@ plus3\p@ minus4\p@
	\textonlyfont@\rm\tenrm
	\textonlyfont@\it\tenit
	\textonlyfont@\sl\tensl
	\textonlyfont@\bf\tenbf
	\textonlyfont@\smc\tensmc
	\ifsyntax@
		\def\big##1{{\hbox{$\left##1\right.$}}}%
	\else
		\let\big\tenbig@
	\textfont0=\tenrm \scriptfont0=\sevenrm \scriptscriptfont0=\fiverm 
	\textfont1=\teni  \scriptfont1=\seveni  \scriptscriptfont1=\fivei
	\textfont2=\tensy \scriptfont2=\sevensy \scriptscriptfont2=\fivesy 
	\textfont3=\tenex \scriptfont3=\tenex \scriptscriptfont3=\tenex
	\textfont\itfam=\tenit \def\it{\fam\itfam\tenit}%
	\textfont\slfam=\tensl \def\sl{\fam\slfam\tensl}%
	\textfont\bffam=\tenbf \def\bf{\fam\bffam\tenbf}%
	\scriptfont\bffam=\sevenbf \scriptscriptfont\bffam=\fivebf
	\textfont\msafam=\tenmsa
	\scriptfont\msafam=\sevenmsa
	\scriptscriptfont\msafam=\fivemsa
	\textfont\msbfam=\tenmsb
	\scriptfont\msbfam=\sevenmsb
	\scriptscriptfont\msbfam=\fivemsb
		\textfont\ttfam=\tentt \def\tt{\fam\ttfam\tentt}%
	\fi
 \normalbaselineskip 14\p@
 \setbox\strutbox=\hbox{\vrule height10\p@ depth4\p@ width0\p@}%
 \normalbaselines\rm \ex@=.2326ex%
}% End of \smalltype.

\def\titletype{%   title fonts
	\def\pointsize@{17}%
	\textonlyfont@\rm\ttlrm
	\ifsyntax@
		\def\big##1{{\hbox{$\left##1\right.$}}}%
	\else
		\let\big\twelvebig@
		\textfont0=\ttlrm \scriptfont0=\twelverm
		\scriptscriptfont0=\tenrm
		\textfont2=\ttlsy \scriptfont2=\twelvesy
		\scriptscriptfont2=\tensy
	\fi
	\normalbaselineskip 25\p@
	\setbox\strutbox=\hbox{\vrule height17\p@ depth8\p@ width0\p@}%
	\normalbaselines
	\rm
	\ex@=.2326ex%
}% End of \titletype.

\def\tenbig@#1{% Same as PLAIN definition of \big .
	{%
		\hbox{%
			$%
			\left
			#1%
			\vbox to8.5\p@{}%
			\right.%
			\n@space
			$%
		}%
	}%
}% End of \tenbig@.

\def\twelvebig@#1{%
	{%
		\hbox{%
			$%
			\left
			#1%
			\vbox to10.2\p@{}%  This is just a guess
			\right.%
			\n@space
			$%
		}%
	}%
}% End of \twelvebig@.

%==============================================================================
%	Macros for symbolic numbering of theorems, figures, etc.
%
%	Purpose:	To provide multiple tracks of symbolic labels for
%			theorems, figures, etc., which will automatically be
%			converted into numbers in a specified format.
%
%	A 'track' or 'type' of labels is identified by a string of letters,
%	with case distinguished.  A symbolic label is a string of any
%	characters.
%
%	Forward references are resolved using two passes and an external file,
%	\jobname.xref.
%
%	Associated with each symbolic label, there are two macros, one holding
%	the actual label, and the other holding the 'state' of the label, which
%	is used for error-checking.
%
%==============================================================================
\newif\ifl@beloutopen
\newwrite\l@belout
\newread\l@belin

\global\let\currentfile=\jobname

% getfile
% =======
% #1: tex file to input
\def\getfile#1{%
	\immediate\closeout\l@belout
	\global\l@beloutopenfalse
	\gdef\currentfile{#1}%
	\input #1%
	\par
	\newpage
}% End of \getfile.

% getxrefs
% ========
% #1: list of other files, separated by commas, from which to grab labels.
\def\getxrefs#1{%
	\bgroup
		\def\gobble##1{}% used to discard \end when we're done
		\edef\list@{#1,}%
		\def\gr@boff##1,##2\end{% process 1st thing in list
			\openin\l@belin=##1.xref
			\ifeof\l@belin
			\else
				\closein\l@belin
				\input ##1.xref
			\fi
			\def\list@{##2}%
			\ifx\list@\empty
				\let\next=\gobble
			\else
				\let\next=\gr@boff
			\fi
			\expandafter\next\list@\end
		}%
		\expandafter\gr@boff\list@\end
	\egroup
}% End of \getxrefs.

% testdefined
% ===========
% #1: control sequence to test
% #2: stuff to execute if #1 is defined
% #3: stuff to execute if #1 is not defined
\def\testdefined#1#2#3{%
	\expandafter\ifx
	\csname #1\endcsname
	\relax
	#3%
	\else #2\fi
}% End of \testdefined

\def\document{%
%	\openout\contents=\jobname.contents
%	\normaltype
	\minaw@11.11128\ex@ % minimum arrow with for @>...>...> macro
%	\pageno=-2 % for preliminary pages; -1 will be for table of cont.
	\def\alloclist@{\empty}%
	\def\fontlist@{\empty}%
	\openin\l@belin=\jobname.xref	% Input \jobname.xref if it exists.
	\ifeof\l@belin\else
		\closein\l@belin
		\input \jobname.xref
	\fi
}% End of \document.

% getst@te
% ========
% #1: type of label
% #2: symbolic label
\def\getst@te#1#2{%
	\edef\st@te{\csname #1s!#2\endcsname}%
	\expandafter\ifx\st@te\relax
		\def\st@te{0}%
	\fi
}% End of \getst@te

% setst@te
% ========
% #1: type of label
% #2: symbolic label
% #3: new state
\def\setst@te#1#2#3{%
	\expandafter
	\gdef\csname #1s!#2\endcsname{#3}%
}% End of setst@te.

% setupautolabel
% ==============
% #1: type of label
% #2: form of actual label (a token list, which should involve the count
%	register \#1Number).
\outer\def\setupautolabel#1#2{%
	\def\newcount@{\global\alloc@0\count\countdef\insc@unt}	% Normally,
		% \newcount is \outer, so we have to repeat the definition
		% here, using some macros from plain.
	\def\newtoks@{\global\alloc@5\toks\toksdef\@cclvi}% See above.
	\expandafter\newcount@\csname #1Number\endcsname
	\expandafter\global\csname #1Number\endcsname=1%
	\expandafter\newtoks@\csname #1l@bel\endcsname
	\expandafter\global\csname #1l@bel\endcsname={#2}%
}% End of \setupautolabel.

% reflabel
% ========
% #1: type of label
% #2: symbolic label
\def\reflabel#1#2{%
	\testdefined{#1l@bel}% See whether the type is known.
	{%	Yes, the type is known.
		\getst@te{#1}{#2}%
		\ifcase\st@te
			%	State 0: the label is undefined.
			???%	Insert dummy label.
			\message{Unresolved forward reference to
				label #2. Use another pass.}%
		\or	%	state 1: defined by file, not yet referenced
			\setst@te{#1}{#2}2%       Why do I need this %?
			\csname #1l!#2\endcsname % Insert the label.
		\or	%	State 2: defined by file, referenced
			\csname #1l!#2\endcsname % Insert the label.
		\or	%	State 3: defined by a \setlabel
			\csname #1l!#2\endcsname % Insert the label.
		\fi
	}{%	No, the type is unknown.
		{\escapechar=-1 % to make the next statement work
		\errmessage{You haven't done a
			\string\\setupautolabel\space for type #1!}%
		}%
	}%
}% End of \reflabel.

{\catcode`\{=12 \catcode`\}=12
	\catcode`\[=1 \catcode`\]=2
	\xdef\Lbrace[{]%	Needed in order to write braces to a file.
	\xdef\Rbrace[}]%
]%

% setlabel
% ========
% #1: type of label
% #2: symbolic label
\def\setlabel#1#2{%
	\testdefined{#1l@bel}%	See whether the type is known.
	{%	Yes, the type is known.
		\edef\templ@bel@{\expandafter\the
			\csname #1l@bel\endcsname}%
		\def\@rgtwo{#2}%
		\ifx\@rgtwo\empty
		\else
			\ifl@beloutopen\else
				\immediate\openout\l@belout=\currentfile.xref
				\global\l@beloutopentrue
			\fi
			\getst@te{#1}{#2}%
			\ifcase\st@te
				% state 0: undefined
			\or	% state 1: defined by file, not yet used
			\or	% state 2: defined by file, referenced
				\edef\oldnumber@{\csname #1l!#2\endcsname}%
				\edef\newnumber@{\templ@bel@}%
				\ifx\newnumber@\oldnumber@
				\else
					\message{A forward reference to label 
						#2 has been resolved
						incorrectly.  Use another
						pass.}%
				\fi
			\or	% state 3: defined by \setlabel
				\errmessage{Same label #2 used in two
					\string\setlabel s!}%
			\fi
			\expandafter\xdef\csname #1l!#2\endcsname
				{\templ@bel@}%	Set the value of the label.
			\setst@te{#1}{#2}3%
			\immediate\write\l@belout % Save label value.
				{\string\expandafter\string\gdef
				\string\csname\space #1l!#2%
				\string\endcsname
				\Lbrace\templ@bel@\Rbrace
				}%
			\immediate\write\l@belout % Save label state.
				{\string\expandafter\string\gdef
				\string\csname\space #1s!#2%
				\string\endcsname
				\Lbrace 1\Rbrace
				}%
		\fi
		\templ@bel@	% Insert the label value.
		\expandafter\ifx\envir@end\endref % inside \ref?
			\gdef\marginalhook@{\marginal{#2}}%
		\else
			\marginal{#2}%    write symbolic label in margin
		\fi
		\expandafter\global\expandafter\advance	% Increment the counter
			\csname #1Number\endcsname
			by 1 %
	}{%	No, the type is unknown.
		{\escapechar=-1
		\errmessage{You haven't done a \string\\setupautolabel\space
			for type #1!}%
		}%
	}%
}% End of \setlabel.

%====================End of symbolic labelling macros==========================
%	Macro for labelling theorems, definitions, etc.
% The following macros allow forward references, at the expense of requiring
% two passes. 

\newcount\SectionNumber
\setupautolabel{t}{\number\SectionNumber.\number\tNumber}
\setupautolabel{r}{\number\rNumber}
\setupautolabel{T}{\number\TNumber}

\define\rref{\reflabel{r}}
\define\tref{\reflabel{t}}

\define\tnum{\setlabel{t}}
\define\rnum{\setlabel{r}}

%%%%%%%%%%%%%%%%%%%%%%%%%%%%%%%%%%%%%%%%%%%%%%%%%%%%%%%%%%%%%%%%%%%%%%%
%  macros for marginal notes
%
\def\strutdepth{\dp\strutbox}%
\def\strutheight{\ht\strutbox}%

\newif\iftagmode
\tagmodefalse

\let\old@tagform@=\tagform@
\def\tagform@{\tagmodetrue\old@tagform@}

\def\marginal#1{%
	\ifvmode
	\else
		\strut
	\fi
	\ifdraftmode
		\ifmmode
			\ifinner
				\let\Vorvadjust=\Vadjust
			\else%				display math mode
				\let\Vorvadjust=\vadjust
			\fi
		\else
			\let\Vorvadjust=\Vadjust
		\fi
		\iftagmode	% special case - tag of an equation
			\llap{%
				\smalltype
				\vtop to 0pt{%
					\pretolerance=2000
					\tolerance=5000
					\raggedright
					\hsize=.72in
					\parindent=0pt
					\strut
					#1%
					\vss
				}%
				\kern.08in
				\iftagsleft@
				\else
					\kern\hsize
				\fi
			}%
		\else% not tagmode
			\Vorvadjust{%
				\kern-\strutdepth % back up to baseline
				{%
					\smalltype
					\kern-\strutheight % match next baseline
					\llap{%
						\vtop to 0pt{%
							\kern0pt
							\pretolerance=2000
							\tolerance=5000
							\raggedright
							\hsize=.5in
							\parindent=0pt
							\strut
							#1%
							\vss
						}%
						\kern.08in
					}%
					\kern\strutheight
				}%
				\kern\strutdepth
			}% end of Vorvadjust.
		\fi% end iftagmode.
	\fi
}% End of \marginal.

% When a \vadjust occurs inside nested boxes, it doesn't seem to do anything,
% so we need to do some trickery to get the adjustment outside.

\newbox\Vadjustbox

\def\Vadjust#1{% My generalization of \vadjust
	\global\setbox\Vadjustbox=\vbox{#1}%
	\ifmmode
		\ifinner
			\innerVadjust
		\fi		%	don't do it in display math mode
	\else
		\innerVadjust
	\fi
}% End of \Vadjust.

\def\innerVadjust{%
	\def\nexti{\aftergroup\innerVadjust}%
	\def\nextii{%
		\ifvmode
			\hrule height 0pt % to prevent \baselineskip glue
			\box\Vadjustbox
		\else
			\vadjust{\box\Vadjustbox}%
		\fi
	}%
	\ifinner
		\let\next=\nexti
	\else
		\let\next=\nextii
	\fi
	\next
}%

\global\let\marginalhook@\empty

\def\endref{%
%      To wind up the preceding box it is convenient to call
%      \makerefbox again; it will also open a new box, however, so we
%      give it arguments \z@ and \endgraf\egroup that will cause the
%      new box to be closed immediately and discarded.
%  \makerefbox\z@{\endgraf\egroup}%
%This is the modification of 3-24-95:
\setbox\tw@\box\thr@@
\makerefbox?\thr@@{\endgraf\egroup}%
%      Then we call \endref@ to take all the saved material and
%      combine it into a paragraph, adding punctuation to separate
%      pieces.
  \endref@
%      The \endgraf is done here rather than in \endref@ because in
%      \moreref or \transl cases \endref@ shouldn't do the \endgraf.
  \endgraf
%      Finally, we need to close the group that was started by \ref.
%      This has the effect of killing the current definition of
%      \envir@end, among other things.
  \endgroup
  \keyhook@
  \marginalhook@
  \global\let\keyhook@\empty % \global to conserve save stack
  \global\let\marginalhook@\empty % \global to conserve save stack
}

\catcode`\@=\oldatsign

%%%%%%%%%%%%%%%%%%%%%%%%%%%%%%%%%%%
%% THE PAPER STARTS HERE  %%%%%%%%%
%%%%%%%%%%%%%%%%%%%%%%%%%%%%%%%%%%%

\nologo
%\draft
\def\fdeg{\operatorname{fdeg}}
\define\Deg{\operatorname{Deg}}
\define\Cl{\Cal C\ell}
\redefine\b{\beta}
\define\e{\varepsilon}
\def \Sym{\operatorname{Sym}}

\define\cok{\operatorname{coker}}
\def \t {\otimes}
\define\({\left(}
\define\){\right)}
\def \p {\oplus}
\def \G{\Bbb G}
\def\F {\Bbb F}
\define\W{\tsize\bigwedge}
\define\tW{\tsize\bigwedge}
\def \rank{\operatorname{rank}}
\def\reg{\operatorname{reg}}
\def\Tor{\operatorname{Tor}}
\def\a{\alpha}
\define \onto {\twoheadrightarrow}
\def \f {\varphi} 
\define\iso{\cong}
\def\depth{\operatorname{depth}}
\def\m{\frak m}
\def\maxm{\frak m}
\def\M{\Bbb M}
\define \g{\gamma}
\def\and{\quad\text{and}\quad}
\def \L {\Bbb L}
\define \Hgy{\operatorname{H}}

\comment

\define\[{\left[}
\define\]{\right]}

\define \maxm {\frak m}
\def\reg{\operatorname{reg}}

\endcomment
\topmatter
\abstract
 Let $A$ be the homogeneous coordinate ring of a rational normal scroll. The ring $A$ is equal to the quotient of a polynomial ring $S$ by the ideal generated by the two by two minors of a scroll matrix $\psi$ with two rows and $\ell$ catalecticant blocks.  
The class group of $A$ is cyclic, and is infinite provided $\ell$ is at least two.   One generator of the class group   is $[J]$, where  $J$ is the ideal of $A$ generated by the entries of the first column of $\psi$. The positive powers of $J$ are well-understood, in the sense that the $n^{\text{th}}$ ordinary power,  the $n^{\text{th}}$ symmetric power, and the $n^{\text{th}}$ symbolic power all coincide and therefore all three $n^{\text{th}}$ powers are resolved by a generalized Eagon-Northcott complex. The inverse of $[J]$ in the class group of $A$ is $[K]$, where $K$ is the ideal generated by the entries of the first row of $\psi$. We study the positive powers of $[K]$. We obtain a minimal generating set and a Gr\"obner basis  for the preimage in $S$ of the symbolic power $K^{(n)}$. We   describe a filtration of   $K^{(n)}$ in which all of the factors   are Cohen-Macaulay $S$-modules resolved by generalized Eagon-Northcott complexes.   We  use this filtration to describe the modules in a  finely graded resolution of $K^{(n)}$ by free $S$-modules. We  calculate the regularity of the graded $S$-module $K^{(n)}$ and we show that the symbolic Rees ring of $K$ is Noetherian.

\endabstract

\title Divisors on Rational Normal Scrolls\endtitle
  \leftheadtext{Kustin, Polini, and Ulrich}
\rightheadtext{Divisors on Rational Normal Scrolls}
 \author Andrew R. Kustin\footnote{Supported in part by the National Security Agency.\hphantom{the National Science Foundation and XXX}}, Claudia Polini\footnote{Supported in part by the National Science Foundation and the National Security Agency.\hphantom{XXX}}, and Bernd  Ulrich\footnote{Supported in part by the National Science Foundation.\hphantom{and the National Security Agency XXX}}\endauthor
 \address
Mathematics Department,
University of South Carolina,
Columbia, SC 29208\endaddress
\email kustin\@math.sc.edu \endemail
\address Mathematics Department,
University of Notre Dame,
Notre Dame, IN 46556\endaddress 
\email     cpolini\@nd.edu \endemail
\address 
Department of Mathematics, 
Purdue University,
West Lafayette, IN 47907
\endaddress
\email          ulrich\@math.purdue.edu \endemail
\endtopmatter

\document

%%%%%%%%%%%%%%%%%%
%% Introduction %%
%%%%%%%%%%%%%%%%%%

\bigpagebreak

\SectionNumber=0\tNumber=1

\flushpar{\bf   Introduction.}\footnote""{\hskip-.17in 2000 {\it Mathematics Subject Classification.} Primary: 13C20, Secondary: 13P10.}\footnote""{\hskip-.17in {\it Key words and phrases.} Divisor Class Group, Filtration, Gr\"obner Basis, Rational Normal Scroll, Regularity, Resolution, Symbolic Rees Ring.}

\medskip

Fix a field $k$ and positive integers $\ell$ and $\sigma_1\ge \sigma_2 \ge \dots\ge \sigma_{\ell}\ge 1$. The rational normal scroll $\operatorname{Scroll}(\sigma_1,\dots,\sigma_{\ell})$ is the image of the map 
$$\pmb \Sigma\:(\Bbb A^2\setminus\{0\})\times (\Bbb A^{\ell}\setminus\{0\})\to \Bbb P^{N},$$where $N=\ell-1+\sum\limits_{i=1}^{\ell}\sigma_i$ and 
$$\pmb \Sigma(x,y;t_1,\dots,t_{\ell})=(x^{\sigma_1}t_1,x^{\sigma_1-1}yt_1,\dots,y^{\sigma_1}t_1,x^{\sigma_2}t_2,x^{\sigma_2-1}yt_2,\dots,y^{\sigma_{\ell}}t_{\ell}).$$ From this     one sees that the homogeneous coordinate ring of $\operatorname{Scroll}(\sigma_1,\dots,\sigma_{\ell}) \subseteq \Bbb P^N$ is the subalgebra 
$$A=k[x^{\sigma_1}t_1,x^{\sigma_1-1}yt_1,\dots,y^{\sigma_1}t_1,x^{\sigma_2}t_2,x^{\sigma_2-1}yt_2,\dots,y^{\sigma_{\ell}}t_{\ell}]$$ of the polynomial ring $k[x,y,t_1,\dots,t_{\ell}]$. This algebra has a presentation  $A=S/I_2(\psi)$, where $S$ is the polynomial ring $$S=k[\{T_{i,j}\mid 1\le i\le \ell\and 1\le j \le \sigma_i+1\}],$$ $\psi$ is the matrix $\psi=\bmatrix \psi_1&\vrule&\dots&\vrule&\psi_{\ell}\endbmatrix$, and for each $u$, $\psi_u$ is 
the generic catalecticant matrix $$\psi_u=\bmatrix T_{u,1}&T_{u,2}&\dots&T_{u,\sigma_u-1}&T_{u,\sigma_u}\\
T_{u,2}&T_{u,3}&\dots&T_{u,\sigma_u}&T_{u,\sigma_u+1}\endbmatrix.$$
Further information about   rational normal scrolls, with alternative descriptions and many applications, may be found in   \cite{\rref{E05},\rref{EH85},\rref{H},\rref{Reid}}.

The class group of the normal domain $A$ is cyclic, and is infinite provided $\ell\ge 2$. One generator of $\Cl(A)$ is $[J]$, where  $J$ is the ideal of $A$ generated by the entries of the first column of $\psi$. The positive powers of $J$ are well-understood, in the sense that the $n^{\text{th}}$ ordinary power $J^n$,  the $n^{\text{th}}$ symmetric power $\Sym_n(I)$, and the $n^{\text{th}}$ symbolic power $J^{(n)}$ all coincide and therefore all three $n^{\text{th}}$ powers are resolved by a  generalized Eagon-Northcott complex. The inverse of $[J]$ in the class group of $A$ is $[K]$, where $K$ is the ideal generated by the entries of the first row of $\psi$. The positive powers of $[K]$ are less well understood. The purpose of the present paper is to rectify this. In Section one we obtain a minimal generating set for $K^{(n)}$; the graded components of this ideal can also be read from \cite{\rref{Sc}, 1.3}. In Section two we exhibit a Gr\"obner basis for the preimage of $K^{(n)}$ in $S$. The Gr\"obner basis is obtained from a    minimal generating set  of $I_2(\psi)$ in $S$ and a monomial minimal generating set of $K^{(n)}$ in $A$. 
In Section three, we   describe a filtration of   $K^{(n)}$ in which all of the factors   are Cohen-Macaulay $S$-modules resolved by generalized Eagon-Northcott complexes.  We  use this filtration in Section four to describe the modules in a  finely graded resolution of $K^{(n)}$ by free $S$-modules. More generally, though less explicitly, resolutions of homogeneous coordinate rings of subvarieties of rational normal scrolls have been approached in \cite{\rref{Sc}, 3.2 and 3.5} in terms of resolutions by locally free sheaves having a filtration by generalized ``Eagon-Northcott sheaves''.   We  calculate the regularity of the graded $S$-module $K^{(n)}$ in Section five. 
The interest in this topic is reflected by the existence of papers like  \cite{\rref{Mi05}} and Hoa's conjecture
\cite{\rref{M00}}.  
In Section six we show that the symbolic Rees ring of $K$ is Noetherian.  (Of course, some symbolic Rees
rings are not Noetherian \cite{\rref{R85},\rref{R90}}, and the question of when the symbolic Rees ring is Noetherian
remains very open \cite{\rref{H87},\rref{C91},\rref{GNW94},\rref{KM08}}.)

Our computation of the regularity in Section five uses a second filtration of $K^{(n)}$ that is coarser than the filtration of Section three. The factors of the coarser filtration are still Cohen-Macaulay modules resolved by generalized Eagon-Northcott complexes. These resolutions give rise to a resolution of $K^{(n)}$ that is sufficiently close to a minimal resolution to allow for a computation of the regularity. On the other hand, if minimality of resolutions is not an issue, like in the calculation  of Hilbert series (see, for example \cite{\rref{kpu-d}}), then it is advantageous to use the finer filtration of Section three as it is easier to describe.

Let $I$ be a homogeneous ideal of height two in $k[x,y]$. Suppose that the presenting matrix of $I$ is almost linear in the sense that the entries of one column have degree $n$ and all of the other entries are linear. In   \cite{\rref{kpu-d}} we prove that the Rees ring and the special fiber ring of $I$ both have the form $A/\Cal A$, where $A$ is the coordinate ring of a rational normal scroll and the ideals $\Cal A$ and $K^{(n)}$ of $A$ are isomorphic. We use the results of the present paper to identify explicit generators for $\Cal A$, to resolve the powers $I^s$ of $I$, to compute the regularity of $I^s$, and to calculate the reduction number of $I$.

%%%%%%%%%%%%%%%% 
%% Section 1  %%
%%%%%%%%%%%%%%%%
%\newpage
\bigpagebreak
\SectionNumber=1\tNumber=1
\flushpar{\bf \number\SectionNumber.\quad  The generators of $K^{(n)}$.}
\medskip
\definition {Data \tnum{:.data37}} We are given integers $\sigma_1\ge \dots\ge \sigma_{\ell}\ge 1$ and an integer $n\ge 2$.
Let $S$ be the polynomial ring $$S=k[\{T_{i,j}\mid 1\le i\le \ell\and 1\le j \le \sigma_i+1\}].$$
For each $u$, with $1\le u\le \ell$, let $\psi_u$ be 
the generic catalecticant matrix $$\psi_u=\bmatrix T_{u,1}&T_{u,2}&\dots&T_{u,\sigma_u-1}&T_{u,\sigma_u}\\
T_{u,2}&T_{u,3}&\dots&T_{u,\sigma_u}&T_{u,\sigma_u+1}\endbmatrix.\tag\tnum{:.psiu}$$  Define $\psi$ to be the matrix
$$\psi=\bmatrix \psi_1&\vrule&\dots&\vrule&\psi_{\ell}\endbmatrix.\tag\tnum{:.psi}$$ 
Let $H$ be the ideal $I_2(\psi)$ of $S$ and $A$  the ring $S/H$. We will write $T_{i,j}$ for a variable in $S$ and also for its image in $A$ -- the meaning will be clear from context. 
Recall that $A$ is a Cohen-Macaulay ring of Krull dimension $\ell+1$ with isolated singularity. In particular, it is a normal domain. 
Let $K$ be the ideal in $A$ generated by the entries of the top row of $\psi$. Notice that $K$ is a height one prime ideal of $A$.
\enddefinition 

In Theorem \tref{T1.1} we identify a generating set for $K^{(n)}$ and in Proposition \tref{L50.3} we identify a minimal generating set for $K^{(n)}$. 

Ultimately, we will put three gradings on the rings $S$ and  $A$. The first grading on $S$ is defined by setting  $$\Deg(T_{i,j})=\sigma_i+1-j.\tag\tnum{Deg}$$  
Notice that $H$ is a homogeneous ideal with respect to this   grading
and thus $\Deg$ induces a grading on $A$, which we also  denote by $\Deg$. Let $A_{\ge n}$ be the ideal of $A$ generated by all monomials $M$ with $\Deg(M)\ge n$.

\proclaim{Theorem \tnum{T1.1}} The $n^{\text{th}}$ symbolic power, $K^{(n)}$, of $K$ is equal to $A_{\ge n}$. \endproclaim
\demo{Proof} Calculate in $A$.  First observe that
$$T_{i,\sigma_{i}+1}^{\sigma_i-j}T_{i,j}=T_{i,\sigma_i}^{\sigma_i+1-j}\in K^{\sigma_i+1-j},\tag\tnum{t1.2}$$ for all $i,j$ with $1\le i\le \ell$ and $1\le j\le \sigma_i$. Indeed, the statement is obvious when $j=\sigma_i$. If $j=\sigma_i-1$, then the assertion holds because 
$$0=\det \bmatrix T_{i,\sigma_i-1}&T_{i,\sigma_i}\\T_{i,\sigma_i}&T_{i,\sigma_i+1}\endbmatrix. $$ The proof of (\tref{t1.2}) is completed by induction on $j$. 
Since $T_{i,\sigma_{i}+1}$ is not in the prime ideal $K$, from
(\tref{t1.2}) we obtain 
$T_{i,j}\in K^{(\sigma_i+1-j)}=K^{(\Deg T_{i,j})}$. Thus,   
$$K^n\subseteq A_{\ge n}\subseteq K^{(n)}.$$

Observe that $\Deg(T_{i,\sigma_i+1})=0$ for $1\le i\le \ell$; hence, $T_{i,\sigma_i+1}$ is regular on $A/A_{\ge n}$ because $A$ is a domain and $n>0$. 
On the other hand, the localization $A[T_{i,\sigma_i+1}^{-1}]$ is a regular ring; and hence, in this ring, $K^{(n)}$ coincides with $K^n$, thus with $A_{\ge n}$.
Since, $T_{i,\sigma_i+1}$ is regular modulo $A_{\ge n}$,
we conclude that $A_{\ge n}$ is equal to  $K^{(n)}$. \qed
\enddemo

\proclaim{Observation \tnum{1.7}}Let $R=k[x,y]$ be a polynomial ring with homogeneous maximal ideal $\maxm$ and write $B=k[x,y,t_1,\dots,t_{\ell}]$. Define the homomorphism of $k$-algebras $\pi\:S\to B$ with $$\pi(T_{i,j})=x^{\sigma_i-j+1}y^{j-1}t_i.$$ 
\roster
\item"{(a)}"The image of $\pi$ is the $k$-subalgebra
$k[R_{\sigma_1}t_1,\dots,R_{\sigma_{\ell}}t_{\ell}]$ of $B$. 
\item"{(b)}"The homomorphism $\pi\:S\to B$ induces an isomorphism $A\iso \pi(S)$.
\item"{(c)}" We have $K\subseteq Bx\cap A=A_{\ge 1}$.
\item"{(d)}" A monomial 
$x^{\alpha}y^{\beta} \prod_{u=1}^{\ell}t_u^{c_u}$ of $B$ belongs to
$\pi(S)$ if and only if $$\alpha + \beta= \sum_{u=1}^{\ell} c_u\sigma_u. \tag\tnum{CLA2}$$
\item"{(e)}" The ring  $\pi(S)$ is a direct summand of $B$ as an $A$-module.
\endroster\endproclaim\demo{Proof}
Assertion (a) is obvious. It follows that the image of $\pi$ is the special fiber ring of the $R$-module $\maxm^{\sigma_1}\p\dots\p\maxm^{\sigma_{\ell}}$. From this one sees that the Krull dimension of $\pi(S)$ is $\ell+1$, which is also the dimension of $A$. Thus, the ideals $H\subseteq \ker \pi$ are prime of the same height and $\pi(S)\iso A$. This is (b).  On $B$, we define a grading by giving $x$ degree $1$ and the other variables degree $0$. The map $\pi$ is homogeneous with respect to this grading on $B$ and the grading $\Deg$ on $S$. Hence, the grading on $B$ induces a grading on the subalgebra $\pi(S)\iso A$ that coincides with $\Deg$ as defined in (\tref{Deg}). Notice that  $K\subseteq Bx\cap A=A_{\ge 1}$, which is (c). (This provides an alternative proof that $K^{(n)}\subseteq Bx^n\cap A=A_{\ge n}$.) 
Assertion (d) is obvious and (e) follows because  the
complementary summand is the $A$-module generated by all monomials of $B$ that
do not satisfy (\tref{CLA2}). \qed
\enddemo

\medskip 
Now we move in the direction of identifying a  minimal generating set for $K^{(n)}$. 
In this discussion, we also use the standard grading, where each variable has degree one, 
we will refer to it as the total degree.
 
\proclaim{Observation \tnum{O50.51}} If $f$ is a monomial of $S$ with $\Deg(f)>0$, then there exists a monomial of the form 
$$M=T_{1,1}^{a_1}\cdots T_{k,1}^{a_k} T_{k,v}T_{k,\sigma_k+1}^{b_k}\cdots T_{\ell,\sigma_{\ell}+1}^{b_{\ell}}\tag\tnum{*1}$$ in $S$ with   $1\le v\le \sigma_k$, 
$\Deg f=\Deg M$, 
and $f-M\in H$.
 \endproclaim
\demo{Proof}We will use this calculation later in the context of Gr\"obner bases; so, we make our argument    very precise. Order the variables of $S$ with
$$T_{1,1}>T_{1,2}>\dots>T_{1,\sigma_1+1}>T_{2,1}>\dots>T_{2,\sigma_2+1}>T_{3,1}>\dots>T_{\ell,\sigma_{\ell}+1}.\tag\tnum{to}$$Observe that there exists $\a\le\b$ such that $f=f_1f'f_2$ where $f_1=T_{1,1}^{a_1}\cdots T_{\a,1}^{a_{\a}}$, $f_2=T_{\b,\sigma_{\b}+1}^{b_{\b}}\cdots T_{\ell,\sigma_{\ell}+1}^{b_{\ell}}$ and $$T_{i,j}|f'\implies T_{\a,1}>T_{i,j}>T_{\b,\sigma_{\b}+1}.$$
Let $T_{i,j}$ be the largest variable which divides $f'$ and $T_{u,v}$ be the smallest variable which divides $f'$. 
We may shrink $f'$, if necessary, and insist that   $1<j$  and $v<\sigma_u+1$. 
If $f'$ has total degree at most one, then one easily may write $f$ in the form of $M$. 
We assume that $f'$ has total degree at least two. Take $f''$ with $f'=T_{i,j}f''T_{u,v}$.
Notice that $$h=-\det \bmatrix T_{i,j-1}&T_{u,v}\\
T_{i,j}&T_{u,v+1}\endbmatrix = T_{i,j}T_{u,v}-T_{i,j-1}T_{u,v+1}\tag\tnum{*2}$$ is in $H$ and 
$$f-f_1hf''f_2=f_1T_{i,j-1}f''T_{u,v+1}f_2\tag\tnum{*3}$$
is more like the desired $M$ than $f$ is. Replacing $f$ by the element of (\tref{*3}) does not change $\Deg$ because the element $h$ of (\tref{*2}) is homogeneous with respect to this grading. Proceed in this manner until $M$ is obtained. \qed \enddemo

We use the notion of eligible tuples when we identify a minimal generating set for $K^{(n)}$ in Proposition \tref{L50.3}. We also use this notion in Section 3 when we describe a filtration of $K^{(n)}$ whose factors are Cohen-Macaulay modules.

\definition{Definition \tnum{D37.8}}\medskip\flushpar {1.} We say that $\pmb a$ is an  {\it  eligible} 
$k$-tuple if $\pmb a$ is a $k$-tuple, $(a_1,\dots,a_k)$, of non-negative integers with $0\le k\le \ell-1$ and   $\sum\limits_{u=1}^{k} a_u\sigma_u<n$.

\medskip\flushpar {2.} Let $\pmb a$ be an eligible $k$-tuple.  The non-negative integer $f(\pmb a)$ is defined by
$$\sum\limits_{u=1}^ka_u\sigma_u+f(\pmb a)\sigma_{k+1} <n\le \sum\limits_{u=1}^ka_u\sigma_u+(f(\pmb a)+1)\sigma_{k+1};$$and the positive integer $r(\pmb a)$ is defined to be
$$r(\pmb a)=\sum\limits_{u=1}^ka_u\sigma_u+(f(\pmb a)+1)\sigma_{k+1}-n+1.$$
Be sure to notice that
$$ 1\le r(\pmb a)\le \sigma_{k+1}.\tag\tnum{:.bstn37}$$
\medskip\flushpar {3.} We write $T^{\pmb a}$ to mean $\prod\limits_{u=1}^kT_{u,1}^{a_u}$ for each eligible $k$-tuple $\pmb a=(a_1,\dots,a_k)$.
\enddefinition
\remark{Remark \tnum{R8}} The empty tuple, $\emptyset$, is always eligible, and we have
$$\tsize f(\emptyset)=\lceil\frac n{\sigma_1}\rceil -1,\quad r(\emptyset)=\sigma_1\lceil\frac n{\sigma_1}\rceil -n+1, \and T^{\emptyset}=1.$$\endremark

\definition{Notation} If $\theta$ is a real number, then $\lceil \theta\rceil$ and $\lfloor \theta \rfloor$ are the ``round up'' and ``round down'' of $\theta$, respectively; that is, $\lceil \theta\rceil$ and $\lfloor \theta \rfloor$ are the integers with
$$ \lceil \theta \rceil-1<\theta \le \lceil \theta \rceil\and
\lfloor \theta \rfloor\le \theta < \lfloor \theta \rfloor+1.
$$\enddefinition

\definition{Definition \tnum{D1.9}}
Let $\Cal L$ be the following list of elements of $S$: $$\Cal L=\bigcup\limits_{k=0}^{\ell-1}\{T^{\pmb a}T_{k+1,1}^{f(\pmb a)}T_{k+1,u}\mid \text{ $\pmb a$ is an eligible $k$-tuple and $1\le u\le r(\pmb a)$}\}.$$\enddefinition

\proclaim{Observation \tnum{O50.5a}}  Let $M$ be the monomial $T_{1,1}^{a_1}\cdots T_{k,1}^{a_k} T_{k,v}T_{k,\sigma_k+1}^{b_k}\cdots T_{\ell,\sigma_{\ell}+1}^{b_{\ell}}$ of $S$. If $\Deg M\ge n$, then   $M$ is divisible by an element of $\Cal L$.  \endproclaim

\demo{Proof}  We have
  $$n\le \Deg(M)=\sum\limits_{u=1}^k a_u\sigma_u+\sigma_k+1-v.$$ 
If $\sum_{u=1}^k a_u\sigma_u<n$, then let $\pmb a$ be the eligible $(k-1)$-tuple $(a_1,\dots,a_{k-1})$. In this case, $f(\pmb a)=a_k$,  $1\le v\le r(\pmb a)$, and  $M$ is divisible by $T^{\pmb a}T_{k,1}^{f(\pmb a)}T_{k,v}
\in \Cal L$. 
 If $n\le \sum_{u=1}^k a_u\sigma_u$, then identify the least index $j$ with $n\le \sum_{u=1}^ja_u\sigma_u$ and let $\pmb a$ be the eligible 
$(j-1)$-tuple $(a_1,\dots,a_{j-1})$. In this case, $f(\pmb a)<a_{j}$
and $M$ is divisible by $T^{\pmb a}T_{j,1}^{f(\pmb a)}T_{j,1}\in \Cal L$.
\qed \enddemo

\proclaim{Observation \tnum{O50.5b}} The ideals $K^{(n)}$ and   $\Cal LA$ are equal. \endproclaim

\demo{Proof} Recall that $K^{(n)}=A_{\ge n}$ according to Theorem \tref{T1.1}. The elements of $\Cal L$ have $\Deg\ge n$, which gives
$\Cal L A\subseteq A_{\ge n}= K^{(n)}$. To prove the other inclusion,  let $f$ be a monomial in $S$ with $\Deg(f)\ge n$. By Observation \tref{O50.51} there exists a monomial $M$  with $f-M\in H$ and $\Deg M=\Deg f \geq n$. Now Observation \tref{O50.5a} shows that $M$ is divisible by an element of $\Cal L$.
\qed \enddemo

\proclaim{Proposition \tnum{L50.3}}
The elements of  $\Cal L$ form  a minimal generating set for
the ideal $K^{(n)}$.
\endproclaim

\demo{Proof}From Observation \tref{O50.5b} we know that $\Cal L$ is a generating set for $K^{(n)}$. To show it is a minimal generating set, we use the map $\pi\: S\to B$ of Observation \tref{1.7} that identifies $A$ with the monomial subring $k[\{x^{\sigma_i-j+1}y^{j-1}t_i\}]$ of $B=k[x,y,t_1,\dots,t_{\ell}]$. The elements of $\pi(\Cal L)$ are monomials in the polynomial ring $B$, and it suffices to show that if $h\in \pi(\Cal L)$ divides $g\in \pi(\Cal L)$ in $B$, then $h=g$ in $B$.

Let $\pmb a$ and $\pmb b$ be eligible $k$ and $j$ tuples, respectively, and let
$$\matrix g=\pi(T^{\pmb a}T_{k+1,1}^{f(\pmb a)}T_{k+1,v})
=x^{G}y^{v-1}t_1^{a_1}\dots t_k^{a_k}t_{k+1}^{f(\pmb a)+1}\qquad \text{and}\hfill \\\vspace{5pt}
h=\pi(T^{\pmb b}T_{j+1,1}^{f(\pmb b)}T_{j+1,w})=
x^{H}y^{w-1}t_1^{b_1}\dots t_j^{b_j}t_{j+1}^{f(\pmb b)+1},\hfill \endmatrix
$$for some $v$ and $w$ with $1\le v\le r(\pmb a)$ and $1\le w\le r(\pmb b)$,
where $$\matrix G=\sum\limits_{u=1}^{k}a_u\sigma_u+(f(\pmb a)+1)\sigma_{k+1}-v+1\qquad \text{and}\hfill\\
H=\sum\limits_{u=1}^{j}b_u\sigma_u+(f(\pmb b)+1)\sigma_{j+1}-w+1.\hfill\endmatrix $$
 The hypothesis that $h$ divides $g$   ensures that $j\le k$ and $b_u\le a_u$ for $1\le u\le j$.   If $j<k$, then $f(\pmb b)+1\le a_{j+1}$ and
$$n\le \sum\limits_{u=1}^{j}b_u\sigma_u+(f(\pmb b)+1)\sigma_{j+1}\le \sum\limits_{u=1}^{j+1}a_u\sigma_u\le \sum\limits_{u=1}^ka_u\sigma_u< n.$$ This contradiction guarantees that $j=k$. Again, the hypothesis ensures that $f(\pmb b)\le f(\pmb a)$, and $b_i\le a_i$, for all $i$. If $b_i<a_i$, for some $i$, then 
$$n\le \sum\limits_{u=1}^{k}b_u\sigma_u+(f(\pmb b)+1)\sigma_{k+1}\le
\sum\limits_{u=1}^{k}a_u\sigma_u+f(\pmb b)\sigma_{k+1}\le
\sum\limits_{u=1}^{k}a_u\sigma_u+f(\pmb a)\sigma_{k+1}<n,$$ since $\sigma_{k+1}\le \sigma_i$. This contradiction guarantees that $\pmb b=\pmb a$.    Again, since $h$ divides $g$, we also have $w\le v$ and $H\le G$. As $\pmb b=\pmb a$, the definition of $H$ and $G$ forces $w=v$. Thus, indeed, $h=g$. 
\qed\enddemo
Inspired by  Observation \tref{1.7} and the proof of Proposition \tref{L50.3}, we introduce the ``fine grading'' on $S$. Let $$\e_u=(0,\dots,0,1,0,\dots,0)\tag\tnum{ms}$$ be the $\ell$-tuple with $1$ in position $u$ and $0$ in all other positions.
The variable $T_{i,j}$ has ``fine degree'' given by
$$\fdeg(T_{i,j})=(\sigma_i-j+1,j-1;\e_i).\tag\tnum{fine}$$ The variables of $S$ have distinct fine degrees. 
Notice that $H$ is homogeneous with respect to fine degree and therefore $\fdeg$ induces a grading on $A$. 
Observe that the grading $\fdeg$ on $A$ 
is simply the grading induced on $A$ by the embedding $A\hookrightarrow B=k[x,y,t_1,\dots,t_{\ell}]$ of Observation \tref{1.7}, where the polynomial ring $B$ is given 
 the usual multigrading.

The two previous gradings that we have considered ($\Deg$ and total degree) can be read from $\fdeg$. Let $\pmb \sigma$ represent the $\ell$-tuple $\pmb \sigma=(\sigma_1,\dots,\sigma_{\ell})$. If $M$ is the monomial 
$$M=\prod\limits_{i=1}^{\ell}\prod\limits_{j=1}^{\sigma_i+1}T_{i,j}^{a_{i,j}}$$ 
of $S$, then 
$$\fdeg M=(\Deg M, \pmb \sigma\cdot \pmb \e-\Deg M; \pmb \e),$$where $\pmb \e$ is the $\ell$-tuple $\pmb \e=( e_1,\dots, e_{\ell})$, with $ e_i=\sum\limits_{j=1}^{\sigma_i+1}a_{i,j}$, and $
\pmb \sigma\cdot \pmb \e$ is the dot product. The total degree of $M$ is $e_1+\dots+e_{\ell}=\pmb 1\cdot \pmb \e$, where $\pmb 1=(1,\dots,1)$ is an $\ell$-tuple of ones. 
We return to the notion of fine degree in (\tref{free}).

%%%%%%%%%%%%%%%% 
%% Section 2  %%
%%%%%%%%%%%%%%%%
%\newpage
\bigpagebreak
\SectionNumber=2\tNumber=1
\flushpar{\bf \number\SectionNumber.\quad  Gr\"obner basis.}
\medskip

In Theorem \tref{T2.1}, we identify a Gr\"obner basis for the preimage of $K^{(n)}$ in $S$; and as an application, in Corollary \tref{C2.1}, we compute $\depth A/K^{(n)}=1$. Sometimes it is convenient to label the variables using a single subscript. That is, we write $T_j$ for $T_{1,j}$; $T_{\sigma_1+1+j}$ for $T_{2,j}$;
$T_{\sigma_1+\sigma_2+2+j}$ for $T_{3,j}$, etc. In this notation, the matrix $\psi$ of (\tref{:.psi}) is
$$\psi=\bmatrix T_1&\dots&T_{\sigma_1}&T_{\sigma_1+2}&\dots&T_{\sigma_1+\sigma_2+1}&T_{\sigma_1+\sigma_2+3}&\dots \\
T_2&\dots&T_{\sigma_1+1}&T_{\sigma_1+3}&\dots&T_{\sigma_1+\sigma_2+2}&T_{\sigma_1+\sigma_2+4}&\dots
\endbmatrix.\tag\tnum{conv}$$ Order the variables of $S$ with $T_1>T_2 >\cdots$, as was done in (\tref{to}).  Impose the reverse lexicographic order on the monomials of $S$. In other words, for two monomials 
$$M_1=T_1^{\a_1}\cdots T_N^{\a_N}\and M_2=T_1^{\b_1}\cdots T_N^{\b_N}$$
one has $M_1>M_2$ if and only if either $\sum\a_i>\sum\b_i$, or else $\sum\a_i=\sum\b_i$ and 
the right most non-zero entry of $(\a_1-\b_1,\dots,\a_N-\b_N)$ is negative. When we study a homogeneous polynomial from $S$ we underline its leading term.  The next result is   well-known, see \cite{\rref{CDR}, Thm. 4.11}. We give a proof for the sake of completeness.  This proof provides  good practice in using the Buchberger criterion for determining when a generating set $G$ of an ideal is a Gr\"obner basis for the ideal. It entails showing that the $S$-polynomial of any two elements of $G$ reduces to zero modulo $G$; see, for example, \cite{\rref{CLO}, Sect\. 2.9, Thm. 3}.

\proclaim{Lemma \tnum{GB1}} The set $G$ of $2\times 2$ minors of $\psi$ forms a Gr\"obner basis for $I_2(\psi)$. \endproclaim
\demo{Proof} Select four columns from $\psi$:
$$\psi'=\bmatrix T_a&T_b&T_c&T_d\\T_{a+1}&T_{b+1}&T_{c+1}&T_{d+1}\endbmatrix,$$with $a\le b\le c\le d$. For $i<j$, let 
$$\Delta_{i,j}=-\det \bmatrix T_i&T_j\\T_{i+1}&T_{j+1}\endbmatrix=\underline{T_{i+1}T_j}-T_iT_{j+1}.$$ 

We first assume that $a<b<c<d$. For most partitions of $\{a,b,c,d\}$ into $p<q$ and $r<s$, the leading terms of $\Delta_{p,q}$ and $\Delta_{r,s}$ are relatively prime; and therefore, the $S$-polynomial $S(\Delta_{p,q},\Delta_{r,s})$   reduces to zero modulo  $G$ (see, for example,  \cite{\rref{CLO}, Sect 2.9, Prop. 4}). The only interesting $S$-polynomial is $S(\Delta_{a,c},\Delta_{b,d})$ when $c=b+1$. In this case, the greatest common divisor of the leading terms of $$\Delta_{a,c}=\underline{T_{a+1}T_c}-T_aT_{c+1}\and \Delta_{b,d}=\underline{T_{b+1}T_d}-T_bT_{d+1}$$ is $T_c=T_{b+1}$; thus
$$S(\Delta_{a,c},\Delta_{b,d})=T_d\Delta_{a,c}-T_{a+1}\Delta_{b,d}=
\underline{-T_aT_{c+1}T_d}+T_{a+1}T_bT_{d+1}.$$We know the generalized Eagon-Northcott complex associated to $\psi'$; and therefore, we know that  the product
$$P=\psi'\bmatrix 0&-\Delta_{c,d}&\Delta_{b,d}&-\Delta_{b,c}\\
\Delta_{c,d}&0&-\Delta_{a,d}&\Delta_{a,c}\\
-\Delta_{b,d}&\Delta_{a,d}&0&-\Delta_{a,b}\\
\Delta_{b,c}&-\Delta_{a,c}&\Delta_{a,b}&0\endbmatrix$$is identically zero.
It follows that 
$$0=-P_{1,2}-P_{2,3}=\left\{\matrix \phantom{-}T_a\Delta_{c,d}-T_c\Delta_{a,d}+T_d\Delta_{a,c}\hfill\\-T_{a+1}\Delta_{b,d}+T_{b+1}\Delta_{a,d}-T_{d+1}\Delta_{a,b},\hfill\endmatrix\right.$$and
$$-T_a\Delta_{c,d}+T_{d+1}\Delta_{a,b}=T_d\Delta_{a,c}-T_{a+1}\Delta_{b,d}=S(\Delta_{a,c},\Delta_{b,d}). \tag\tnum{t2.8'}
$$
The leading term of each summand of the left hand side of (\tref{t2.8'}) is at most the leading term of the right hand side; hence,  the $S$-polynomial $S(\Delta_{a,c},\Delta_{b,d})$ reduces to zero modulo  $G$.

There are no complicated calculations to make if some of the indices $a,b,c,d$ are equal. Indeed, it suffices to consider these cases:
$$\matrix
a=b<c<d&\implies&S(\Delta_{a,c},\Delta_{a,d})=T_d\Delta_{a,c}-T_c\Delta_{a,d}=-T_a\Delta_{c,d}\hfill\\
a<b=c<d&\implies&\text{the leading terms of $\Delta_{a,b}$ and $\Delta_{b,d}$ are relatively prime} \hfill\\
a<b<c=d&\implies&S(\Delta_{a,c},\Delta_{b,c})=T_{b+1}\Delta_{a,c}-T_{a+1}\Delta_{b,c}=T_{c+1}\Delta_{a,b}.\hfill\endmatrix$$In each case, the relevant $S$-polynomial reduces to zero modulo  $G$. \qed \enddemo
Retain the notation of (\tref{:.data37}). Recall the polynomials $G$ from   Lemma \tref{GB1} and $\Cal L$ from Definition \tref{D1.9}.

\proclaim{Theorem \tnum{T2.1}}The set of polynomials $G\cup \Cal L$ in $S$ is a Gr\"obner basis for the preimage of $K^{(n)}$ in $S$. \endproclaim

\demo{Proof} Again we apply the Buchberger criterion. We saw in Lemma \tref{GB1} that every $S$-polynomial $S(h_1,h_2)$, with $h_1,h_2\in G$, reduces to zero modulo  $G\cup \Cal L$. If $M_1,M_2$ are in $\Cal L$, then the $S$-polynomial $S(M_1,M_2)$ is equal to zero. Finally, we study the $S$-polynomial $f=S(M_1,h_1)$, where $M_1$ is an element of $\Cal L$ and $h_1$ is in $G$. The only interesting case is when $M_1$ and the leading term of $h_1$ have a factor in common. Henceforth, we make this assumption. It is clear that $f$ is monomial. We claim that $\Deg(f)\ge n$. Once the claim is established, then  Observation  \tref{O2.4} shows that $f$ reduces to zero modulo  $G\cup \Cal L$. We prove the claim.
Write 
$$h_1=-\det\bmatrix T_{i,j-1}&T_{u,v}\\T_{i,j}&T_{u,v+1}\endbmatrix=\underline{T_{i,j}T_{u,v}}-T_{i,j-1}T_{u,v+1}$$ for variables $T_{i,j-1}>T_{i,j}\ge T_{u,v}>T_{u,v+1}$ from $S$. 
There are three possibilities for the greatest common divisor of $M_1$ and $T_{i,j}T_{u,v}$:
$$ T_{i,j}\quad\text{or}\quad T_{u,v}\quad\text{or}\quad T_{i,j}T_{u,v}.$$
In the first case,  $f=\frac{M_{1}}{T_{i,j}}T_{i,j-1}T_{u,v+1}$ and
$$\Deg(f)-\Deg(M_{1})=1+\Deg(T_{u,v+1})\ge 1.$$  In the second case, $f=\frac{M_{1}}{T_{u,v}}T_{i,j-1}T_{u,v+1}$ 
and
$$\Deg(f)-\Deg(M_{1})=-1+\Deg(T_{i,j-1})\ge 0.$$ In the third case,   
$f=\frac{M_{1}}{T_{i,j}T_{u,v}}T_{i,j-1}T_{u,v+1}$ and $\Deg(f)=\Deg(M_{1})$. 
In each case, $\Deg(f)\ge \Deg(M_1)\ge n$. Thus the claim is established and the proof is complete.
\qed
\enddemo
 \proclaim{Observation  \tnum{O2.4}} If $f$ is a monomial of $S$ with $\Deg(f)\ge n$, then $f$ reduces to zero modulo  $G\cup \Cal L$.\endproclaim
\demo{Proof} The proof of Observation \tref{O50.51} shows that the remainder of $f$ on division by $G$ has the form of $M$ from (\tref{*1}) with $\Deg(M)=\Deg(f) \geq n$. (The proof of Observation \tref{O50.51} does not mention division by $G$; however, the binomial $h$ of (\tref{*2}) is in $G$ and the leading term of    $h$   is $T_{i,j}T_{u,v}$. This leading term divides the only term of $f$ with quotient
$$\frac{f}{T_{i,j}T_{u,v}}=f_1f''f_2.$$ We
calculate the $S$-polynomial $S(f,h)=f-hf_1f''f_2$ in (\tref{*3}). Proceed in this manner until $M$ is obtained.) 
Furthermore, Observation \tref{O50.5a} shows that $M$ is divisible by an element of $\Cal L$. \qed \enddemo

\proclaim{Corollary \tnum{C2.1}} Adopt the notation of Data \tref{:.data37} with $n\ge 2$, then $\depth A/K^{(n)}=1$. \endproclaim

\demo{Proof}Let $\kappa$ represent the preimage of $K^{(n)}$ in $S$. We compute $\depth S/\kappa$. The variable $T_{\ell,\sigma_{\ell}+1}$ does not divide the leading term of any element  in the Gr\"obner basis $G\cup \Cal L$ of $\kappa$.  Therefore $T_{\ell,\sigma_{\ell}+1}$ is regular on $S/\kappa$.

Now we show that the homogeneous maximal ideal $\maxm$ of $S$ is an associated prime of
$S/(\kappa,T_{\ell,\sigma_{\ell}+1})=C$. Indeed, $\maxm$ annihilates the image of the element  $T_{\ell,\sigma_{\ell}}$ in $C$, because the images of $K, T_{\ell,\sigma_{\ell}+1}, H$ are all zero in $C$. On the other hand,  $T_{\ell,\sigma_{\ell}}$ maps to a nonzero element in $C=A/ (A_{\geq n}, T_{\ell,\sigma_{\ell}+1}A)$. 
For if $T_{\ell,\sigma_{\ell}}A \subseteq (A_{\geq n}, T_{\ell,\sigma_{\ell}+1}A)$, then
$T_{\ell,\sigma_{\ell}}A\subseteq T_{\ell,\sigma_{\ell}+1}A$ because $\Deg T_{\ell,\sigma_{\ell}}=1$ and $A_{\geq n}$ is generated by homogeneous elements with $\Deg \ge n>1$. But  $T_{\ell,\sigma_{\ell}}A \subset T_{\ell,\sigma_{\ell}+1}A$ is impossible since $A=S/H$ with $H \subset \maxm^2$.
\qed \enddemo

%%%%%%%%%%%%%%%% 
%% Section 3  %%
%%%%%%%%%%%%%%%%
%\newpage
\bigpagebreak
\SectionNumber=3\tNumber=1
\flushpar{\bf \number\SectionNumber.\quad  Filtration.}
\medskip

In Theorem \tref{:.O37.11}, we   describe a filtration of the $n^{\text{th}}$ symbolic power, $K^{(n)}$, of $K$. The factors in this filtration   are Cohen-Macaulay $S$-modules.   We  use this filtration to describe the modules in a  $\fdeg$-graded resolution of $K^{(n)}$ by free $S$-modules, see Theorem \tref{:.Tmay2} and (\tref{fine}). We  calculate the regularity of the graded $S$-module $K^{(n)}$ in Theorem \tref{:Reg2}. 

\definition{Definition \tnum{:.D37.8}}Recall the notation of Definition \tref{D37.8}.
 
\medskip\flushpar {1.}   We put a total order on the set of  eligible tuples. If $\pmb b=(b_1,\dots,b_j)$ and $\pmb a=(a_1,\dots,a_k)$ are eligible  tuples, then we say that $\pmb b>\pmb a$ if either
$$\matrix \text{(a)\  $j<k$ and $b_i=a_i$ for $1\le i\le j$, or}\hfill\\
\text{(b)\   $\exists i\le \min\{j,k\}$ with $b_i>a_i$ and  $b_s=a_s$ for $1\le s\le i-1$.}\hfill\endmatrix\tag\tnum{:.>}$$  
 If one pretends that $\pmb b$ and $\pmb a$ have the same length, filled out as necessary on the right by the symbol $\infty$: $(b_1,\dots,b_j,\infty,\dots,\infty)$ and $(a_1,\dots,a_k,\infty,\dots,\infty)$, then one may test the  total order $>$ of (\tref{:.>}) using only rule (b).  Recall from Remark \tref{R8} that the empty tuple $\emptyset$ is always an eligible tuple. Notice that $\emptyset$ is the largest eligible tuple.
 
\medskip\flushpar {2.} For an eligible tuple $\pmb a$ we define the $A$-ideals 
$$\split &\Cal D_{\pmb a}=\sum_{\pmb b>\pmb a} T^{\pmb b}T_{j+1,1}^{f(\pmb b)}(T_{j+1,1},\dots,T_{j+1,r(\pmb b)}) \quad\text{and}\\&\Cal E_{\pmb a}=\sum_{\pmb b\ge \pmb a} T^{\pmb b}T_{j+1,1}^{f(\pmb b)}(T_{j+1,1},\dots,T_{j+1,r(\pmb b)}),\endsplit$$where $\pmb b=(b_1,\dots, b_j)$ is  eligible and $j$ is arbitrary. Notice that $\Cal D_{\emptyset}=0$, and if the tuple $\pmb a$ is not empty, then  
   $$ \Cal D_{\pmb a}=\sum_{\pmb b >\pmb a} \Cal E_{\pmb b},$$ where the sum is taken over all eligible tuples $\pmb b$ with $\pmb b>\pmb a$. Notice also that if 
$\pmb a$ is an eligible $k$-tuple, then $$\Cal E_{\pmb a}=\Cal D_{\pmb a}+T^{\pmb a}T_{k+1,1}^{f(\pmb a)}(T_{k+1,1},\dots,T_{k+1,r(\pmb a)}).$$ 
This gives a finite filtration 
$$(0)\subsetneq \Cal E_{\emptyset}\subsetneq \dots \subsetneq \Cal E_{0^{\ell-1}}=K^{(n)},$$ of $K^{(n)}$, where $0^s$ is the $s$-tuple $(0,\dots,0)$.
We define two parallel collections of ideals $\{\Cal E_{\pmb a}\}$ and $\{\Cal D_{\pmb a}\}$ simultaneously because there is no convenient way to denote the eligible tuple which is immediately larger than a particular eligible tuple $\pmb a$. Notice that the modules $\Cal E_{\pmb a}/\Cal D_{\pmb a}$ are exactly the factors of the filtration $\{\Cal E_{\pmb a}\}$. 
\enddefinition

 Recall the fine grading (\tref{fine}) on $S$ and $A$. Observe that the ideals $\Cal D_{\pmb a}$ and $\Cal E_{\pmb a}$ are homogeneous in this grading. Define   $\fdeg$-graded free $S$-modules
   $$E=\left\{\matrix S(1,-1;0)\\\p\\ S(0,0;0)\endmatrix \right.\qquad\and\qquad F_u=\left\{\matrix 
S(-\sigma_u+1,-1;-\e_u)\\\p\\
S(-\sigma_u+2,-2;-\e_u)\\\p\\ \vdots\\\p\\ 
S(0,-\sigma_u;-\e_u),\endmatrix \right.\tag\tnum{free}$$ for $1\le u\le \ell$. Notice that each $\psi_u\:F_u\to E$ is a homogeneous map, with respect to $\fdeg$; and therefore, for each $k$, the cokernel
of 
$$\psi_{>k}=\bmatrix \psi_{k+1}&\mid&\dots&\mid&\psi_{\ell}\endbmatrix\:\bigoplus\limits_{u=k+1}^{\ell}F_u\to E\tag\tnum{next}$$ is a graded $S$-module, with respect to the $\fdeg$-grading.
Let $\pmb a$ be an eligible $k$-tuple. In this section we prove that $\Cal E_{\pmb a}/\Cal D_{\pmb a}$ is a well-known Cohen-Macaulay module. Let   $P_k$ be the ideal 
$$P_{k}=(\{T_{i,j}\mid 1\le i\le k\and 1\le j\le \sigma_i+1\})$$ of $S$, and $\e_{\pmb a}$ be the multi-shift
$$\e_{\pmb a}=\sum\limits_{u=1}^ka_u\e_u,$$for $\e_u$ given in (\tref{ms}).
In Theorem \tref{:.O37.11} we prove that the $\fdeg$-graded $S$-modules $\Cal E_{\pmb a}/\Cal D_{\pmb a}$ and  
$$\Sym_{r(\pmb a)-1}^{S/I_2(\psi_{>k})}(\cok(\psi_{>k}))\t_S S/P_{k}\(-\pmb \sigma\cdot \pmb \e,0;-\pmb\e\)\tag\tnum{:.ED}$$ 
are isomorphic, for $\pmb \e=\e_{\pmb a}+(f(\pmb a)+1)\e_{k+1}$.
The module of (\tref{:.ED}) might look more familiar if we observe that 
$$ \matrix\Sym_{r(\pmb a)-1}^{S/I_2(\psi_{>k})}(\cok(\psi_{>k}))\t_S \frac S{P_{k}}\hfill \\\vspace{5pt} \cong (T_{k+1,1},T_{k+1,2})^{r(\pmb a)-1}\frac A{P_kA}((r(\pmb a)-1)(\sigma_{k+1},0;\e_{k+1}));\hfill \endmatrix $$ see the proof of Lemma \tref{:L42.L}.
We have written $\Sym_{r(\pmb a)-1}^{S/I_2(\psi_{>k})}$ rather than $\Sym$ or $\Sym^S$ in order to emphasize that when $r(\pmb a)-1=0$, then the module of (\tref{:.ED}) is a shift of 
$$S/I_2(\psi_{>k})\t_S S/P_{k}=A/P_kA.$$
Recall, from (\tref{:.bstn37}), that $r(\pmb a)-1$  is non-negative and is  less than the number of columns of $\psi_{>k}$. 
The ideal $I_2(\psi_{>k})$ has generic height (equal to the number of columns of $\psi_{>k}$ minus $1$) and the symmetric power $r(\pmb a)-1$ is small enough that  $\Sym_{r(\pmb a)-1}(\cok(\psi_{>k}))$ is a perfect $S$-module and is resolved by a generalized Eagon-Northcott complex. (See, for example, the family of complexes studied in and near Theorem A2.10 in \cite{\rref{E95}} or Theorem 2.16 in \cite{\rref{BV88}}. Recall that the $S$-module $M$ is perfect if the grade of the annihilator of $M$ on $S$ is equal to the projective dimension of $M$.)

The module (\tref{:.ED}) is annihilated by $P_k$. The first step in the proof of Theorem \tref{:.O37.11} is to show that $\Cal E_{\pmb a}/\Cal D_{\pmb a}$ is also annihilated by $P_k$.

\proclaim{Lemma \tnum{:+.L37.!}} If $\pmb a$ is an eligible $k$-tuple and
$(\kappa, r)$ is a pair of  integers with $$k+1\le \kappa\le \ell\and
1\le r\le \sum\limits_{u=1}^{k}a_u\sigma_u+f(\pmb a)\sigma_{k+1}-n+1+\sigma_{\kappa},\tag\tnum{ub}$$
then $P_kT^{\pmb a}T_{k+1,1}^{f(\pmb a)}T_{\kappa,r}\subseteq  \Cal D_{\pmb a}$. In particular, if $\pmb a$ is an eligible $k$-tuple, then
\roster\item"{\rm(a)}"$P_k\Cal E_{\pmb a}
\subseteq \Cal D_{\pmb a}$, and
\item"{\rm(b)}" if
 $r(\pmb a)=\sigma_{k+1}$, then
 $P_kT^{\pmb a}T_{k+1,1}^{f(\pmb a)}T_{\kappa,r}\subseteq  \Cal D_{\pmb a}$, for all $(\kappa,r)$ with $k+1\le \kappa\le \ell$ and   $1\le r\le \sigma_{\kappa}$.\endroster
\endproclaim
\demo{Proof}
We notice that (a) and (b) are applications of the first assertion.
Indeed, if $\kappa=k+1$, then the upper bound on $r$ in (\tref{ub}) is equal to $r(\pmb a)$;
furthermore, $$\Cal E_{\pmb a}=
T^{\pmb a}T_{k+1,1}^{f(\pmb a)}(\{T_{k+1,r}\mid 1\le r\le r(\pmb a)\})+\Cal D_{\pmb a}.$$
In (b), the hypothesis $r(\pmb a)=\sigma_{k+1}$ forces $\sum\limits_{u=1}^{k}a_u\sigma_u+f(\pmb a)\sigma_{k+1}=n-1$,
and in this case the bound on $r$ in (\tref{ub}) becomes $1\le r\le \sigma_{\kappa}$.

We prove the first assertion.
Fix   $i$ and $s$  with $1\le i\le k$ and   $1\le s\le \sigma_i+1$.
 Let $$X=T_{i,s}T^{\pmb a}T_{k+1,1}^{f(\pmb a)}T_{\kappa,r}.$$
We will prove that $X\in \Cal D_{\pmb a}$.

Define $a_{k+1}=f(\pmb a)$ and
$$b_u=\cases a_u&\text{if $1\le u\le k+1$ and $u\neq i$, and }\\ a_{i}+1&\text{if $u=i$.}\endcases$$
Notice that for each $u$, with $1\le u\le k$, we have
$$(b_1,\dots,b_{u})>\pmb a,$$where we define order as in Definition \tref{:.D37.8}.1. We know
$$\sum\limits_{u=1}^{k}a_u\sigma_u+f(\pmb a)\sigma_{k+1}<n\le \sum\limits_{u=1}^{k}a_u\sigma_u+f(\pmb a)\sigma_{k+1}+\sigma_{k+1}\le \sum\limits_{u=1}^{k+1}b_u\sigma_u.$$Select the least integer $j$ with
$$n\le \sum\limits_{u=1}^{j}b_u\sigma_u.$$ Notice that $i\le j\le k+1$. Select the largest value $b_j'$ with
$$\sum\limits_{u=1}^{j-1}b_u\sigma_u+b_j'\sigma_j< n.$$Notice that
$$0\le b_j'<b_j.$$

Let $\pmb b=(b_1,\dots,b_{j-1})$.
We see that $\pmb b$ is an eligible $(j-1)$-tuple and $\pmb b>\pmb a$.
We have chosen $b_j'$ so that $b_j'=f(\pmb b)$. It follows that
$$T^{\pmb b}T_{j,1}^{b_j'}(T_{j,1},\dots,T_{j,r(\pmb b)})\subseteq \Cal E_{\pmb b}\subseteq \Cal D_{\pmb a}.$$
Write $\rho=\min\{r(\pmb b), s+r-1\}.$ Since $1 \leq \rho
\leq r(\pmb b)$, it suffices to prove that
$$X\in T^{\pmb b}T_{j,1}^{b_j'}T_{j,\rho}A.\tag\tnum{CLA}$$
Notice that

\roster
\item"{(\tnum{a})\phantom{3}}" if $i<j$, then $b_j'<b_j=a_j$, and

\item"{(\tnum{b})}" if $i=j$, then $b_j'=a_j$.
\endroster

 \smallskip We prove (\tref{b}). The definition of $b_j'$ says that
$b_j'$ is the largest integer with
$$\sum_{u=1}^{j-1}b_u\sigma_u+b_j'\sigma_i<n.$$
In other words, $b_j'$ is the largest integer with
$\sum_{u=1}^{j-1}a_u\sigma_u+b_j'\sigma_j<n$. On the other hand, we
know
$$\sum_{u=1}^{j-1}a_u\sigma_u+a_j\sigma_j=\sum_{u=1}^ja_u\sigma_u<n\le \sum_{u=1}^jb_u\sigma_u= \sum_{u=1}^{j-1}a_u\sigma_u+(a_j+1)\sigma_j.$$The last equality holds because $i=j$; so, $b_j=b_i=a_i+1=a_j+1$.  Assertion (\tref{b})
is established.

To prove (\tref{CLA}) we use the embedding $A\hookrightarrow
B=k[x,y,t_1, \ldots, t_{\ell}]$ induced by the map $\pi$ of Observation
\tref{1.7}. Thus (\tref{CLA}) is equivalent to showing that
$$
x^{\g}y^{s+r-2}t_it_{\kappa}\prod_{u=1}^{k+1}t_u^{a_u}
= F x^{\delta}y^{\rho-1}t_j^{b_j'+1}\prod_{u=1}^{j-1}t_u^{b_u},$$
for some $F\in A$, with
$$\g= \sum_{u=1}^{k+1}a_u\sigma_u + \sigma_i+
\sigma_{\kappa}-s-r+2\and \delta = \sum_{u=1}^{j-1}b_u\sigma_u+(b_j'+1)\sigma_j
-\rho+1.$$
Clearly such $F$ exists in the quotient field
of $A$. According to (\tref{a}) and (\tref{b}) one has
$$F=\cases x^{\alpha}y^{\beta}t_{\kappa}\prod\limits_{u=j+1}^{k+1}t_u^{a_u}&\text{if $i=j$}\\x^{\alpha}y^{\beta}t_{\kappa}t_j^{a_j-b_j'-1}\prod\limits_{u=j+1}^{k+1}t_u^{a_u}&\text{if $i<j$,}\endcases$$
and $F$ is an element of $k(x,y)[t_1, \ldots, t_{\ell}]$.
Notice that
$$\beta=s+r-\rho-1\and
\alpha=\cases \nu
&\text{if $i=j$}\\
\nu+\sigma_j(a_j-b_j'-1)&\text{if
$i<j$,}\endcases$$  for $$\nu=\sum_{u=j+1}^{k+1}a_u\sigma_u +
\sigma_{\kappa}-s-r+\rho+1.$$

Recall that $F$ is in the quotient field of $A$ and that $A$ is a direct 
summand of $B$ according to Observation 1.7. Thus, to prove that
$F \in A$ is suffices to show that $F \in B$ or, equivalently,
$$\alpha
\geq 0 \and  \beta \geq 0.\tag\tnum{CLA5}$$
Clearly, $\beta \geq 0$ by the definition of $\rho$. Likewise, if $\rho=s+r-1$, then
$\alpha \geq 0$ according to (\tref{a}). Thus we
may assume that $\rho=r(\pmb b)$. Use the definition of $r(\pmb b)$:
$$r(\pmb b)=\sum\limits_{u=1}^{j-1}b_u\sigma_u+(b_j'+1)\sigma_j-n+1.$$ Treat the cases $i=j$ and $i<j$ separately. Two straightforward calculations yield
$$\a= 
\(\sum_{u=1}^{k+1}a_u\sigma_u+\sigma_{\kappa}-n + 1 -r \) +\(
\sigma_i+1-s \)\geq 0,$$ where the first summand is non-negative by assumption (\tref{ub}) and the second summand is non-negative because of the choice of $s$.   This completes the proof of
(\tref{CLA5}). \qed
\enddemo

We have established   half of Theorem  \tref{:.O37.11}. The next two Lemmas are used in the other 
half of the proof. 
\proclaim{Lemma \tnum{:.O37.9}} If $\pmb a$ is an eligible $k$-tuple, then 
$\Cal D_{\pmb a}\subseteq P_kA$ and 
$T^{\pmb a}$  is not zero in $(A/\Cal D_{\pmb a})_{P_kA}$.
 \endproclaim \demo{Proof} 
Let  
 $\frak S$ be the multiplicative subset of  $A\setminus P_kA$ which consists of the non-zero elements of the ring 
$$\frac {k[T_{k+1,*},\dots,T_{\ell,*}]}{I_2\(\psi_{>k}\)}\tag\tnum{:.kappa}$$ and let $Q$ be the quotient field of the ring of (\tref{:.kappa}). 
We notice that $$\frak S^{-1}(A)=\frac {Q[T_{1,*},\dots,T_{k,*}]}{HQ[T_{1,*},\dots,T_{k,*}]}.$$ Furthermore, since $k\le \ell-1$, $HQ[T_{1,*},\dots,T_{k,*}]$ is generated by linear forms, and $T_{i,j}$ is an associate of $T_{i,1}$ in $\frak S^{-1}(A)$, for all $i,j$ with $1\le i\le k$ and $1\le j \le \sigma_i+1$. Indeed, $\frak S^{-1}(A)$ is naturally isomorphic  to the polynomial ring $Q[T_{1,1},\dots,T_{k,1}]$ in $k$ variables over the field $Q$. Observe that the ring $A_{P_kA}$ is equal to the further localization $Q[T_{1,1},\dots,T_{k,1}]_{(T_{1,1},\dots,T_{k,1})}$ of $\frak S^{-1}(A)$.

We first  show that $T^{\pmb a}$ is not zero in $\frak S^{-1}(A/\Cal D_{\pmb a})$. We have seen that $\frak S^{-1}(A)$ is the polynomial ring $Q[T_{1,1},\dots,T_{k,1}]$. We now observe that  the ideal $\Cal D_{\pmb a}$ of $\frak S^{-1}(A)$ is generated by the following set of monomials:
$$ \matrix \{T^{\pmb b}T_{j+1,1}^{f(\pmb b)+1}\mid \text{ $\pmb b=(b_1,\dots,b_j)$ is eligible, $j<k$, and (\tref{:.>}.a) or (\tref{:.>}.b) is in effect}\}\\ {}\cup\{ T^{\pmb b}\mid  \text{ $\pmb b=(b_1,\dots,b_k)$ is eligible  and 
(\tref{:.>}.b) is in effect}\}.\endmatrix$$It is obvious that none of the monomials in the second set can divide $T^{\pmb a}$. If some monomial from the first set divides $T^{\pmb a}$, then the definition of $f(\pmb b)$, together with the fact that $\pmb a$ is eligible, yields: 
$$n\le \sum_{u=1}^jb_u\sigma_u+(f(\pmb b)+1)\sigma_{j+1}\le \sum_{u=1}^ka_u\sigma_u<n,$$ and of course, this is impossible. 
Thus, $T^{\pmb a}$ is not zero in $\frak S^{-1}(A/\Cal D_{\pmb a})$, which is 
a standard graded $Q$-algebra. We localize at the homogeneous maximal ideal  to see that $T^{\pmb a}$ is also not zero in 
$(A/\Cal D_{\pmb a})_{P_kA}$. In particular, this is not the zero ring, showing that $\Cal D_{\pmb a}\subseteq P_kA$. 
\qed \enddemo

 \proclaim{Lemma \tnum{:L42.L}} Let $\pmb a$ be an eligible $k$-tuple, $B$  the ring $A/\Cal D_{\pmb a}$, and $J$  the ideal $T_{k+1,1}^{f(\pmb a)}(T_{k+1,1},\dots,T_{k+1,r(\pmb a)})$ of $A$. Then \roster 
\item"{\rm(a)}" the module of {\rm (\tref{:.ED})} is isomorphic to $J\frac A{P_kA}(-\pmb \sigma\cdot\e_{\pmb a},0; -\e_{\pmb a})$, and
\item"{\rm(b)}"  $\Cal E_{\pmb a}/\Cal D_{\pmb a}=T^{\pmb a}JB$.\endroster
\endproclaim
\demo{Proof}  Assertion (b) is clear. We prove (a) by establishing  the following sequence     of isomorphisms:
$$\matrix \tsize \Sym_{r(\pmb a)-1}(\cok(\psi_{>k}))\t_S \frac S{P_{k}}\hfill\\\vspace{3pt}\tsize\phantom{XXX}@>\a_1>> (T_{k+1,1},T_{k+1,2})^{r(\pmb a)-1}\frac A{P_kA}((r(\pmb a)-1)(\sigma_{k+1},0;\e_{k+1}))\hfill\\\vspace{3pt}\tsize \phantom{XXX}@>\a_2>> (T_{k+1,1},\dots,T_{k+1,r(\pmb a)})\frac A{P_kA}(\sigma_{k+1},0;\e_{k+1})\hfill\\\vspace{3pt}\tsize \phantom{XXX}@>\a_3>> J\frac A{P_kA}((f(\pmb a)+1)(\sigma_{k+1},0;\e_{k+1})).\hfill\endmatrix $$
The ideal $(T_{k+1,1},T_{k+1,2})$ of the domain $\frac A{P_kA}$ is generated by the entries of the first column of $\psi_{>k}$. The map
$$E@> \bmatrix T_{k+1,2}&-T_{k+1,1}\endbmatrix>> A(\sigma_{k+1},0,\e_{k+1}),$$for $E$ in (\tref{free}), induces  a natural surjection $$\tsize \cok(\psi_{>k})\t_S \frac S{P_{k}}\onto (T_{k+1,1},T_{k+1,2})\frac A{P_kA}(\sigma_{k+1},0;\e_{k+1}).\tag\tnum{:.ns}$$
The map $\a_1$ is the surjection induced by (\tref{:.ns}).
Recall that the domain of $\a_1$ is Cohen-Macaulay and that it has rank one as a module over the domain $A/P_kA$. Furthermore, the target of $\a_1$  is, up to shift, a non-zero ideal in this domain. It follows that $\a_1$ is an isomorphism.
The ideals 
$$(T_{k+1,1},T_{k+1,2})^{r(\pmb a)-1} \and T_{k+1,1}^{r(\pmb a)-2}(T_{k+1,1},\dots,T_{k+1,r(\pmb a)})$$ of the domain $\frac A{P_kA}$ are equal; and therefore, the isomorphism $\a_2$ is given by multiplication by the unit $1/T_{k+1,1}^{r(\pmb a)-2}$ in the quotient field of $\frac A{P_kA}$.
Multiplication by the non-zero element $T_{k+1,1}^{f(\pmb a)}$ of the domain $\frac A{P_kA}$ gives the $A$-module isomorphism $\a_3$. \qed\enddemo

The next lemma is the  final step in our proof of Theorem \tref{:.O37.11}. We will also use the same lemma in the proof of Proposition \tref{:.O37.1*}.

\proclaim{Lemma \tnum{:L44.70}}
Let $\pmb a$ be an eligible $k$-tuple, $B$   the ring $A/\Cal D_{\pmb a}$,    and $\Cal J$   an ideal  of $A$.
Assume that 
\itemitem{\rm(1)} $\Cal J$ is $\fdeg$-homogeneous and the generators of $\Cal J$ involve only the variables $\{T_{i,j}\}$ with $i\ge k+1$, and  
\itemitem{\rm(2)} $P_k$ annihilates $T^{\pmb a}\Cal JB$. 

\flushpar Then  the graded $A$-modules
$T^{\pmb a}\Cal JB$ and $\Cal J(A/P_kA)  (-\pmb\sigma\cdot\e_{\pmb a},0;-\e_{\pmb a})$ are isomorphic.\endproclaim 
\demo{Proof} We exhibit $A$-module homomorphisms
$$\tsize \a\: T^{\pmb a}\Cal JB(\pmb\sigma\cdot\e_{\pmb a},0;\e_{\pmb a})\to \Cal J(\frac A{P_kA})\and \b\: \Cal J(\frac A{P_kA})\to T^{\pmb a}\Cal JB(\pmb\sigma\cdot\e_{\pmb a},0;\e_{\pmb a}),$$
which are inverses of one another.

We first show that $\b\: \Cal J(\frac A{P_kA})\to T^{\pmb a}\Cal JB(\pmb\sigma\cdot\e_{\pmb a},0;\e_{\pmb a})$, given by $\b(X)=T^{\pmb a}X$, for all $X$ in $\Cal J$, is a well-defined $A$-module homomorphism. 
Consider the composition 
$$A\to B\to T^{\pmb a}B(\pmb\sigma\cdot\e_{\pmb a},0;\e_{\pmb a}),$$where the first map is the natural quotient map and the second map is multiplication by $T^{\pmb a}$. This composition restricts to give $\b'\:\Cal J A\to T^{\pmb a}\Cal JB(\pmb\sigma\cdot\e_{\pmb a},0;\e_{\pmb a})$. 
The first hypothesis ensures that $\Cal J A\cap P_kA=\Cal JP_k A$  and the second hypothesis ensures that
$\Cal JP_k A\subseteq \ker \b'$.  So, $\b'$ induces
$$\b\: \Cal J\(\frac A{P_kA}\)= \frac{\Cal J A}{\Cal JA\cap P_kA}\to  T^{\pmb a}\Cal JB(\pmb\sigma\cdot\e_{\pmb a},0;\e_{\pmb a}),$$as described above.

Now we show that
$\a\: T^{\pmb a}\Cal JB(\pmb\sigma\cdot\e_{\pmb a},0;\e_{\pmb a})\to \Cal J(\frac A{P_kA})$, given by  $\a(T^{\pmb a}X)=X$, for all $X$ in $\Cal J$,
is a  well-defined $A$-module homomorphism. Let $$\f\:B=\frac A{\Cal D_{\pmb a}}\to \frac A{P_kA}$$ be the natural quotient map which is induced by the inclusion $\Cal D_{\pmb a}\subseteq P_kA$ of Lemma \tref{:.O37.9} and let $\pi\: B\to T^{\pmb a}B(\pmb\sigma\cdot\e_{\pmb a},0;\e_{\pmb a})$ be multiplication by $T^{\pmb a}$.

The kernel of $\pi$ is the annihilator of $T^{\pmb a}$ in $B$, and the kernel of $\f$ is $P_kB$. We saw in Lemma \tref{:.O37.9} that $T^{\pmb a}\neq 0$ in $B_{P_kB}$. It follows that 
$$\ker \pi\subseteq \ker \f. $$ 
Thus, there exists   a unique $A$-module homomorphism $\f'\: T^{\pmb a}B(\pmb\sigma\cdot\e_{\pmb a},0;\e_{\pmb a})\to \frac A{P_kA}$ for which the diagram 
$$\xymatrix{ B  \ar[r]^{\f}\ar[d]^{\pi}&
 \frac A{P_kA}\\  
T^{\pmb a}B(\pmb\sigma\cdot\e_{\pmb a},0;\e_{\pmb a})\ar[ur]^{\f'} }$$
commutes. The restriction of $\f'$ to $T^{\pmb a}\Cal JB(\pmb\sigma\cdot\e_{\pmb a},0;\e_{\pmb a})$ is the homomorphism $\a$ which is described above. \qed\enddemo

 The next result follows from Lemma \tref{:L42.L} and Lemma \tref{:L44.70}, applied to the ideal $T_{k+1,1}^{f(\pmb a)}(T_{k+1,1},\dots,T_{k+1,r(\pmb a)})$ of $A$; notice that assumption (2) of Lemma \tref{:L44.70} is satisfied according to Lemma \tref{:+.L37.!}(a).
\proclaim{Theorem \tnum{:.O37.11}} Adopt the hypotheses of 
\tref{:.data37}.
  Let $\{\Cal E_{\pmb a}\}$, 
as $\pmb a$ varies over all eligible tuples, be the filtration of $K^{(n)}$ from Definition \tref{:.D37.8}.
 Then, for each   eligible $k$-tuple $\pmb a$,  
the
$\fdeg$-graded $S$-modules $\Cal E_{\pmb a}/\Cal D_{\pmb a}$ and {\rm (\tref{:.ED})} are isomorphic.  \endproclaim

%%%%%%%%%%%%%%%% 
%% Section 4  %%
%%%%%%%%%%%%%%%%
%\newpage
\bigpagebreak
\SectionNumber=4\tNumber=1
\flushpar{\bf \number\SectionNumber.\quad Resolution.}
\medskip

We first record the minimal homogeneous resolution of the module $\Cal E_{\pmb a}/\Cal D_{\pmb a}$ by free $\fdeg$-graded $S$-modules. Recall the free $\fdeg$-graded $S$-modules $E$ and $F_u$ of (\tref{free}). These modules have rank $2$ and $\sigma_u$, respectively. Let  
 $F=F_1\p\dots\p F_\ell$.
The matrices $\psi$ and $\psi_u$ of (\tref{:.psi}) and (\tref{:.psiu}) describe homogeneous $\fdeg$-graded homomorphisms $\psi\:F\to E$ and  $\psi_u\:F_u\to E$. Let  $G_u$ be the free $\fdeg$-graded  $S$-module$$G_u=\left\{\matrix
S(-\sigma_u,0;-\e_u)\\\p\\S(-\sigma_u+1,-1;-\e_u)\\\p\\\vdots \\\p\\S(0,-\sigma_u;-\e_u)\endmatrix\right.$$ of rank $\sigma_{u+1}$,  and let
   ${\rho_u\:G_u\to S}$ be the $\fdeg$-homogeneous $S$-module homomorphism given by 
  $$\rho_u=\bmatrix T_{u,1}&T_{u,2}&\dots&T_{u,\sigma_u+1}\endbmatrix.$$ 
For any $k$ with $0\le k\le \ell-1$, let $F_{>k}$ and $G_{\le k}$ be the free $\fdeg$-graded $S$-modules
$$F_{>k}= \bigoplus\limits_{u=k+1}^\ell F_u \and  G_{\le k}=\bigoplus\limits_{u=1}^kG_u,$$
and let $\psi_{>k}\: F_{>k}\to E$ and $\rho_{\le k}\:G_{\le k}\to S$ be the $\fdeg$-homogeneous $S$-module homomorphisms
$$ \psi_{>k}=\bmatrix \psi_{k+1}&\vrule&\dots&\vrule&\psi_{\ell}\endbmatrix\and  \rho_{\le k}=\bmatrix \rho_{1}&\vrule&\dots&\vrule&\rho_{k}\endbmatrix.$$
  The Koszul complex
$$\G_{k,\bullet}={\tsize\W^{\bullet}G_{\le k}},$$ associated to $\rho_{\le k}\:G_{\le k}\to S$, is a homogeneous resolution of  $S/P_{k}$ by free $\fdeg$-graded $S$-modules. 
We see that 
$$\G_{k,q}=\sum\limits_{i_1+\dots+i_k=q}\tW^{i_1}G_1\t\dots\t\tW^{i_k}G_k\quad \text{for $0\le q\le \sum\limits_{i=1}^k(\sigma_i+1).$}$$
The generalized Eagon-Northcott complex $\F_{\pmb a,\bullet}$, where
$$\F_{\pmb a,p}=\cases \Sym_{r(\pmb a)-1-p}E\t \tW^pF_{>k}&\text{if $0\le p\le r(\pmb a)-1$}\\
D_{p-r(\pmb a)}E^*\t \tW^{p+1}F_{>k}&\text{if $r(\pmb a)\le p\le \rank F_{>k}-1$},\endcases
$$ is a homogeneous resolution of  $\Sym_{r(\pmb a)-1}^{S/I_2(\psi_{>k})}(\cok(\psi_{>k}))$ 
 by free $\fdeg$-graded $S$-modules. See, for example,   \cite{\rref{E95}, Theorem A2.10} or \cite{\rref{BV88}, Theorem 2.16}.
One other generalized Eagon-Northcott complex is of interest to us. For each integer $k$, with $0\le k\le \ell-1$, the complex $(\F_{k,\bullet},d_{k,\bullet})$, with
$$\F_{k,p}=D_{p}E^*\t \tW^{p+1}F_{>k},$$ 
is a homogeneous resolution of 
$$\tsize (\{T_{\kappa,r}\mid k+1\le \kappa\le \ell\ \text{and}\ 1\le r\le \sigma_{\kappa}\})\frac S{I_2(\psi_{>k})}(1,-1;0) \tag \tnum{shft}$$ by $\fdeg$-graded free $S$-modules. The complex $\F_{k,\bullet}$ is called $\Cal C^{-1}$ in \cite{\rref{E95}}. The $\fdeg$-homogeneous augmentation map from the complex $\F_{k,\bullet}$ to the module of
(\tref{shft}) 
is induced by the map
$$\F_{k,0}=F_{>k} @> \bmatrix \xi_{k+1}&\xi_{k+2}&\dots &\xi_{\ell}\endbmatrix >> \frac S{I_2(\psi_{>k})}(1,-1;0),$$  
where $\xi_u\: F_u\to S(1,-1;0)$ is the $\fdeg$-homogeneous map given by  
$$\bmatrix T_{u,1}&T_{u,2}&\dots&T_{u,\sigma_{u}} \endbmatrix$$ and the free $\fdeg$-graded  module $F_u$ is described in (\tref{free}).

 With respect to total degree,
the maps in $\F_{\pmb a,\bullet}$ are linear everywhere, except $\F_{\pmb a,r(\pmb a)}\to \F_{\pmb a,r(\pmb a)-1}
$, where the maps are quadratic because they involve $2\times 2$ minors of $\psi_{>k}$. All of the maps in $\F_{k,\bullet}$ are linear. In other words, with respect to total degree, 
$$\reg \Sym_{r(\pmb a)-1}^{S/I_2(\psi_{>k})}(\cok(\psi_{>k}))=\cases 0&\text{if $k=\ell-1$ and $r(\pmb a)=\sigma_{\ell}$}\\1&\text{in all other cases}. \endcases$$ (A thorough discussion of regularity may be found in Section 5.) Furthermore,
$$\reg \tsize (\{T_{\kappa,r}\mid k+1\le \kappa\le \ell\ \text{and}\ 1\le r\le \sigma_{\kappa}\})\frac S{I_2(\psi_{>k})}=1$$ because the generators live in degree one and the resolution is linear.

\proclaim{Observation \tnum{:.O28'}}Let $k$ be an integer with $0\le k\le \ell-1$.
\roster\item"{\rm(a)}"If $\pmb a$ is an eligible $k$-tuple,  then 
$$(\Bbb L_{\pmb a,\bullet},d_{\pmb a,\bullet})=\(\F_{\pmb a,\bullet}\t_S \G_{k,\bullet}\)\(
-\pmb \sigma\cdot \pmb \e,0;-\pmb \e\)$$
is 
the minimal homogeneous $\fdeg$-graded resolution of the module of $\Cal E_{\pmb a}/\Cal D_{\pmb a}$
by free $S$-modules, for $\pmb \e=\e_{\pmb a}+(f(\pmb a)+1)\e_{k+1}$. 
\item"{\rm(b)}" The complex 
$\Bbb L_{k,\bullet}=\F_{k,\bullet}\t_S \G_{k,\bullet} $ is the minimal homogeneous $\fdeg$-graded resolution of the module  
$$\tsize (\{T_{\kappa,r}\mid k+1\le \kappa\le \ell\ \text{and}\ 1\le r\le \sigma_{\kappa}\})\frac A{P_kA}(1,-1;0)\tag\tnum{lst} $$ by free $S$-modules. 
\item"{\rm(c)}" The $S$-module  $\Cal E_{\pmb a}/\Cal D_{\pmb a}$ 
and  the $S$-module  of {\rm (\tref{lst})}
are Cohen-Macaulay and perfect  of projective dimension $\sum\limits_{u=1}^{\ell}\sigma_u+k-1$.
\endroster\endproclaim

\demo{Proof} Recall that
$$\Cal E_{\pmb a}/\Cal D_{\pmb a}\iso
\Sym_{r(\pmb a)-1}^{S/I_2(\psi_{>k})}(\cok(\psi_{>k}))\t_S S/P_{k}\(-\pmb \sigma\cdot \pmb \e,0;-\pmb\e\)$$ 
by Theorem \tref{:.O37.11}.
 We know that $\F_{\pmb a,\bullet}$ is a minimal homogeneous $\fdeg$-graded resolution of $\Sym_{r(\pmb a)-1}^{S/I_2(\psi_{>k})}(\cok(\psi_{>k}))$ and $\G_{k,\bullet}$ is a resolution of $S/P_{k}$. Furthermore, the generators of $P_k$ are a regular sequence on the  $S$-module $\Sym_{r(\pmb a)-1}^{S/I_2(\psi_{>k})}(\cok(\psi_{>k}))$; therefore,
$$\Tor_i^S(\Sym_{r(\pmb a)-1}^{S/I_2(\psi_{>k})}(\cok(\psi_{>k})), S/P_{k})=0\text{ for all $i\ge 1$,}$$ and
$\F_{\pmb a,\bullet}\t_S \G_{k,\bullet}$ is a minimal homogeneous $\fdeg$-graded resolution of $$\Sym_{r(\pmb a)-1}^{S/I_2(\psi_{>k})}(\cok(\psi_{>k}))\t_SS/P_{k}.$$ 
Notice that the length of this resolution is
$$\sum_{u=k+1}^{\ell}\sigma_u-1 +\sum_{u=1}^k(\sigma_u+1)=\sum_{u=1}^{\ell}\sigma_u+k-1,$$ which is the grade  of the annihilator of the module it resolves. Assertion (a) and half of assertion (c) have been established. The rest of the result 
  is proved the same manner. \qed
\enddemo

Finally, we resolve $K^{(n)}$. Let
$$M=M_0\supseteq M_1 \supseteq \dots \supseteq M_s=0$$ 
be
 a filtration of a module $M$.  If one can resolve each sub-quotient $M_i/M_{i+1}$, then one can resolve $M$ by  an iterated application of the Horseshoe Lemma, as explained in 
 Lemma \tref{:.IMC'}. 
We apply the lemma to the filtration $\{\Cal E_{\pmb a}\}$ of  $K^{(n)}$ in Theorem \tref{:.Tmay2}. One may also apply the lemma to the  filtration $\{\Cal E'_{\pmb a}\}$ of Section 5 without any difficulty.
Neither resolution is minimal.
\proclaim{Lemma \tnum{:.IMC'}} Let $M$ be a finitely generated multi-graded module over a multi-graded Noetherian ring and let $$M=M_0\supseteq M_1 \supseteq \dots \supseteq M_s=0$$ 
be
 a finite filtration by graded submodules. 
Suppose that for each $i$, with $0\le i\le s-1$,
$$\F_{i,\bullet}:\quad \cdots @>d_{i,2} >>\F_{i,1}@>d_{i,1} >>\F_{i,0}$$ is a homogeneous resolution of $M_i/M_{i+1}$. Then, for each $i,j,k$, with $0\le i\le s-1$, $1\le k\le s-i-1$, and $1\le j$, there exists a  homogeneous map 
$$\a_{i,j}^{(k)}\:\F_{i,j}\to \F_{i+k,j-1}$$ such that 
$$(\M,D)\:\quad \dots\to \M_2@>D_2>>\M_1@>D_1>>\M_0
$$ is a homogeneous resolution of $M$, where $\M_j=\bigoplus\limits_{i=0}^{s-1} \F_{i,j}$ and $D_j\:\M_j\to \M_{j-1}$ is the lower triangular matrix
$$D_j=\bmatrix 
d_{0,j}&0&0&0&\dots&0&0\\
\a_{0,j}^{(1)}&d_{1,j}&0&0&\dots&0&0\\
\a_{0,j}^{(2)}&\a_{1,j}^{(1)}&d_{2,j}&0&\dots&0&0\\
\a_{0,j}^{(3)}&\a_{1,j}^{(2)}&\a_{2,j}^{(1)}&d_{3,j}&\dots&0&0\\
\vdots& \vdots&\vdots&\vdots&&\vdots&\vdots\\
\a_{0,j}^{(s-1)}&\a_{1,j}^{(s-2)}&\a_{2,j}^{(s-3)}&\a_{3,j}^{(s-4)}&\dots&\a_{s-2,j}^{(1)}&d_{s-1,j}\endbmatrix.$$
\endproclaim
\demo{Proof} By iteration, it suffices to treat the case $s=2$. 
In this case the proof is a graded version of the Horseshoe Lemma.  \qed\enddemo

\proclaim{Theorem \tnum{:.Tmay2}}Adopt the hypotheses of \tref{:.data37} and recall the resolution $(\Bbb L_{\pmb a,\bullet},d_{\pmb a,\bullet})$ of Observation \tref{:.O28'}. For each triple $(\pmb a,\pmb b,j)$, where  $j$ is a positive integer and $\pmb b>\pmb a$ are eligible tuples,  there exists an $\fdeg$-homogeneous  $S$-module homomorphism 
$$\a_{\pmb a,\pmb b,j}\:\Bbb L_{\pmb a,j}\to \Bbb L_{\pmb b,j-1},$$
such that $$(\frak L,D)\:\quad 0\to \frak L_s\to\dots\to \frak L_2@>D_2>>\frak L_1@>D_1>>\frak L_0
$$ is an $\fdeg$-homogeneous  resolution of $K^{(n)}$, where  ${s=\sum\limits_{u=1}^{\ell}\sigma_u+\ell-2}$, $\frak L_j$ is equal to $\bigoplus\limits_{\pmb a} \Bbb L_{\pmb a,j}$, and  the component 
$$\Bbb  L_{\pmb a,j}\hookrightarrow \frak L_{j}@> D_j >> \frak L_{j-1}@>\operatorname{proj} >> \Bbb L_{\pmb c,j-1}$$
of the map $D_j\:\frak L_j\to \frak L_{j-1}$ is  
is equal to 
$$\cases 
0&\text{if $\pmb a>\pmb c$}\\
d_{\pmb a,j} &\text{if $\pmb a=\pmb c$}\\
\a_{\pmb a,\pmb c,j}&\text{if $\pmb c>\pmb a$}. \endcases$$
\endproclaim
\demo{Proof} Apply Lemma \tref{:.IMC'} to the filtration $\{\Cal E_{\pmb a}\}$ of $K^{(n)}$. \qed\enddemo
\remark{Remark} The length of the complex $\L_{0^{\ell-1},\bullet}$ is $\sum\limits_{u=1}^{\ell}\sigma_u+\ell-2$, which is the same as the projective dimension of $K^{(n)}$ as an $S$-module, as may be calculated from the Auslander-Buchsbaum formula. Indeed, Corollary \tref{C2.1} shows that the depth of $K^{(n)}$, as an $S$-module, is $2$ and it is clear that $S$ has depth equal to $\sum\limits_{u=1}^{\ell}\sigma_u+\ell$.
\endremark
%%%%%%%%%%%%%%%% 
%% Section 5  %%
%%%%%%%%%%%%%%%%
%\newpage
\bigpagebreak
\SectionNumber=5\tNumber=1
\flushpar{\bf \number\SectionNumber.\quad Regularity.}
\medskip
  
We turn our attention to the Castelnuovo-Mumford
regularity of $K^{(n)}$. In this discussion all of the variables of the polynomial ring $S$ have degree one. In Section one, we referred to this situation as the grading on $S$ is given by ``total degree''.  If 
$M$ is a finitely generated non-zero graded $S$-module and 
$$0\to F_k\to \dots \to F_0\to M\to 0,$$ with $F_{i}=\bigoplus_{j=1}^{t_i}S(-a_{i,j})$, is the minimal homogeneous resolution of $M$ by free $S$-modules, then 
the regularity  of $M$  is equal to  $$\reg(M) =\max_{i,j} \{a_{i,j} - i\}
=\max\{n\mid \Hgy_{\maxm}^i(M)_{n-i}\neq 0\text{ for some $i\ge 0$}\},
$$ where $\maxm$ is the maximal homogeneous ideal of $S$. For $M=0$ one sets $\reg(M)=-\infty$.

 There are two contributions to the regularity of $K^{(n)}$. 
The highest generator degrees of $K^{(n)}$ and of $\Cal E_{0^{\ell-1}}/\Cal D_{0^{\ell-1}}$ coincide, where $0^{\ell-1}$ is the $(\ell-1)$-tuple of zeros. Also, most of the generalized Eagon-Northcott complexes $\F_{\pmb a,\bullet}$ are linear in all positions except one position where the maps are quadratic. The rest of the generalized Eagon-Northcott complexes are linear in all positions. For example, the generators of 
$\Cal E_{0^{\ell-1}}/\Cal D_{0^{\ell-1}}$ have degree $\lceil \frac n{\sigma_{\ell}}\rceil$ and the complex $\L_{0^{\ell-1},\bullet}$ contains some quadratic maps if and only if $\sigma_{\ell}\not|(n-1)$. It follows that 
$$\reg\(\Cal E_{0^{\ell-1}}/\Cal D_{0^{\ell-1}}\)=\left.\cases \lceil \frac n{\sigma_{\ell}}\rceil+1&\text{if $\sigma_{\ell}\not|n-1$}\\
  \lceil \frac n{\sigma_{\ell}}\rceil&\text{if $\sigma_{\ell}|n-1$}\endcases\right\} = \left\lceil \frac{n-1}{\sigma_{\ell}}\right\rceil +1.\tag\tnum{new}$$ 
We prove in Theorem \tref{:Reg2} that $\reg K^{(n)}=\reg\(\Cal E_{0^{\ell-1}}/\Cal D_{0^{\ell-1}}\)$. The filtration $\{\Cal E_{\pmb a}\}$ is too fine to allow us to read the exact value of $\reg K^{(n)}$ directly from the factors of  $\{\Cal E_{\pmb a}\}$. In order to complete our calculation of $\reg K^{(n)}$, we introduce a second filtration $\{\Cal E'_{\pmb a}\}$, with $\{\Cal E_{\pmb a}\}$ a refinement of $\{\Cal E'_{\pmb a}\}$.  
 
\definition{Definition \tnum{:D42.66}} The $k$-tuple $\pmb a$ is eligible$'$ if $\pmb a$ is eligible and either $k=\ell-1$ or $r(\pmb a)<\sigma_{k+1}$. If $\pmb a$ is an eligible$'$ $k$-tuple, then 
\roster\item"{(1)}" $\Cal E'_{\pmb a}=\Cal E_{\pmb a}$, and
\item"{(2)}" $\Cal D'_{\pmb a}=\sum \Cal E'_{\pmb b}$, where the sum varies over all eligible$'$ tuples $\pmb b$,  with $\pmb b>\pmb a$.\endroster\enddefinition

Notice that the modules $\Cal E'_{\pmb a}/\Cal D'_{\pmb a}$ are exactly the factors of the filtration $\{\Cal E'_{\pmb a}\}$. The next result, about the filtration $\{\Cal E'_{\pmb a}\}$,  is comparable to Theorem \tref{:.O37.11} about the filtration $\{\Cal E_{\pmb a}\}$. 
From the point of view of regularity, Proposition \tref{:.O37.1*} says that the factors 
$\Cal E'_{\pmb a}/\Cal D'_{\pmb a}$ of the filtration $\{\Cal E'_{\pmb a}\}$ are either factors 
$\Cal E_{\pmb a}/\Cal D_{\pmb a}$ of the filtration $\{\Cal E_{\pmb a}\}$ or else have linear resolution. We delay the proof of Proposition \tref{:.O37.1*} until after we have used the result to prove Theorem \tref{:Reg2}. 

\proclaim{Proposition \tnum{:.O37.1*}}    Let $\pmb a$ be an eligible$'$ $k$-tuple. The $S$-module $\Cal E'_{\pmb a}/\Cal D'_{\pmb a}$ is Cohen-Macaulay and perfect.
\itemitem{\rm (a)} If $r(\pmb a)<\sigma_{k+1}$, then $\Cal E'_{\pmb a}/\Cal D'_{\pmb a}=\Cal E_{\pmb a}/\Cal D_{\pmb a}$ and the assertions of Theorem \tref{:.O37.11}  apply.
\itemitem{\rm (b)} If $r(\pmb a)=\sigma_{k+1}$, then there exists a non-negative integer $j$  such that 
there is an isomorphism of  $\fdeg$-graded $S$-modules: $$\Cal E'_{\pmb a}/\Cal D'_{\pmb a}\iso J(A/P_jA)(-\pmb \sigma\cdot \pmb \e,0;-\pmb \e),$$
where  
$J$ is the $A$-ideal generated by the entries in the first row of $\psi_{>j}$ and $\pmb \e=\e_{\pmb a}+f(\pmb a)\e_{k+1}$.  Furthermore, the complex $\L_{j,\bullet}(-\pmb \sigma\cdot \pmb \e-1,1;-\pmb \e)$ of  Observation \tref{:.O28'} is a resolution of 
$\Cal E'_{\pmb a}/\Cal D'_{\pmb a}$.  
 If all of the variables of $S$ are given degree one, then the minimal $S$-resolution of $\Cal E'_{\pmb a}/\Cal D'_{\pmb a}$ is linear. 
\itemitem{\rm (c)} The modules $\Cal E'_{0^{\ell-1}}/\Cal D'_{0^{\ell-1}}$ and $\Cal E_{0^{\ell-1}}/\Cal D_{0^{\ell-1}}$ are equal.
\endproclaim

 \proclaim{Lemma \tnum{:Reg1}}
Let $R$ be a standard graded Noetherian ring over a field, $M$ a non-zero
finitely generated graded $R$-module, and  $M=M_0 \supseteq M_1 \supseteq
\ldots \supseteq M_s=0$ a finite filtration by graded modules with
factors $N_i=M_i/M_{i+1}$. If ${\operatorname {reg}}\, N_0 \geq {\operatorname {reg}} \,
N_i$ for every $i$, then ${\operatorname {reg}} \, M={\operatorname {reg}} \, N_0$ and ${\depth}\, M \leq {\dim} \, N_0$.
\endproclaim
\demo{Proof} Notice that ${\operatorname {reg}} \, M_{1} \leq {\operatorname {reg}} \, N_0$ and
${\operatorname {reg}} \, M \leq {\operatorname {reg}} \, N_0$. Let $\m$ be the maximal
homogeneous ideal of $R$ and let $d$ be such that
$[\Hgy_{\m}^{d}(N_0)]_{{\operatorname {reg}} N_0 -d} \not= 0$. Clearly $ d \leq
{\dim} \, N_0$. We claim that $[\Hgy_{\m}^{d}(M)]_{{\operatorname {reg}} N_0 -d}
\not= 0$, which gives ${\operatorname {reg}} \,  M \geq {\operatorname {reg}} \, N_0  $ as
well as ${\depth}\, M \leq d \leq {\dim} \, N_0$.

Suppose $[\Hgy_{\m}^{d}(M)]_{{\operatorname {reg}} N_0 -d} = 0$, then the short
exact sequence
$$
0 \longrightarrow M_1 \longrightarrow M \longrightarrow N_0
\longrightarrow 0
$$
induces an embedding
$$
0 \not=[\Hgy_{\m}^{d}(N_0)]_{{\operatorname {reg}} N_0 -d} \hookrightarrow
[\Hgy_{\m}^{d+1}(M_1)]_{{\operatorname {reg}} N_0 -d}.
$$
Hence  $[\Hgy_{\m}^{d+1}(M_1)]_{{\operatorname {reg}} N_0 -d} \not= 0$, which gives
$ {\operatorname {reg}} \, M_1 \geq {\operatorname {reg}} \, N_0 +1$.
\qed
\enddemo

\proclaim{Theorem \tnum{:Reg2}}
Adopt the hypotheses of \tref{:.data37}. Then ${\operatorname {reg}} \, K^{(n)}=\left\lceil\frac{n-1}{\sigma_{\ell}}
\right\rceil + 1$.
\endproclaim
\demo{Proof} Consider the finite
filtration $\{\Cal E'_{\pmb a}\}$ of $K^{(n)}$ as described in Definition \tref{:D42.66}. The factors of this filtration are denoted $\Cal E'_{\pmb a}/\Cal D'_{\pmb a}$, as $\pmb a$ varies over all eligible$'$ tuples. Notice that $0^{\ell-1}$ is the smallest eligible$'$-tuple and $K^{(n)}/ \Cal D'_{0^{\ell -1}}= \Cal E'_{0^{\ell -1}}/ \Cal D'_{0^{\ell -1}}$
has   regularity
$\left\lceil\frac{n-1}{\sigma_{\ell}} \right\rceil + 1$ by Proposition \tref{:.O37.1*}(c) and (\tref{new}). Hence by
Lemma~\tref{:Reg1} it suffices to show that ${\operatorname {reg}} \, \Cal E'_{\pmb a}/\Cal D'_{\pmb a}
\leq \left\lceil\frac{n-1}{\sigma_{\ell}} \right\rceil + 1$ for
every eligible$'$ $k$-tuple $\pmb a$.

The module $\Cal E'_{\pmb a}/\Cal D'_{\pmb a}$ is generated in degree $\sum a_i +
f(\pmb a)+1$ and hence, according to Proposition \tref{:.O37.1*}, has regularity equal to
$$\cases \sum a_i +f(\pmb a)+2,&\text{if $r(\pmb a)<\sigma_{k+1}$,}\\
\sum a_i +f(\pmb a)+1,&\text{if $r(\pmb a)=\sigma_{k+1}$.}\endcases$$ If $r(\pmb a)<\sigma_{k+1}$, then 
$\sum a_i\sigma_i+f(a)\sigma_{k+1}<n-1$.  The hypothesis 
 $\sigma_1 \geq \ldots \geq \sigma_{\ell}$ ensures that
$\sum a_i +f(a) <\frac{n-1}{\sigma_{\ell}}$, and hence $\reg  (\Cal E'_{\pmb a}/\Cal D'_{\pmb a})\le \left\lceil\frac{n-1}{\sigma_{\ell}} \right\rceil + 1$. On the other hand, if 
$r(\pmb a)=\sigma_{k+1}$, then $\sum a_i +f(a) \le\frac{n-1}{\sigma_{\ell}}$, and we still have ${\operatorname {reg}} \( \Cal E'_{\pmb a}/\Cal D'_{\pmb a}\)
\le \left\lceil\frac{n-1}{\sigma_{\ell}} \right\rceil + 1$.
 \qed
\enddemo

\remark{Remark}Recall that $\Cal E_{0^{\ell-1}}'/\Cal D_{0^{\ell-1}}'=\Cal E_{0^{\ell-1}}/\Cal D_{0^{\ell-1}}$ according to Proposition \tref{:.O37.1*}(c). Since the later module has dimension two, Lemma \tref{:Reg1} and the proof of Theorem \tref{:Reg2} yield an alternate proof of Corollary \tref{C2.1}: $\depth K^{(n)}=2$ for $n\ge 2$. \endremark

\medskip
We begin our proof of Proposition \tref{:.O37.1*} by making a more detailed study of the totally ordered set of all eligible tuples. In particular, sometimes it is clear when a pair of  eligible tuples are adjacent.   

\definition {Notation \tnum{:N42}} If $\pmb a$ is an eligible $k$-tuple with $k<\ell-1$, then let $N(\pmb a)$ be the  $(k+1)$-tuple $(\pmb a,f(\pmb a))$. If $2\le h< \ell-k$, then let $N^h(\pmb a)=N(N^{h-1}(\pmb a))$. We let $N^0$ denote the identity function.
\enddefinition

\proclaim{Lemma \tnum{:sd42}} Let $\pmb a$ be an eligible $k$-tuple with $k<\ell-1$.
\roster
\item"{\rm (a)}" The $(k+1)$-tuple $N(\pmb a)$ is eligible and  the eligible tuples $\pmb a>N(\pmb a)$ are nearest  neighbors in the sense that if $\pmb b$ is an eligible tuple with $\pmb a\ge \pmb b\ge N(\pmb a)$, then either $\pmb a=\pmb b$ or $\pmb b =N(\pmb a)$.
\item"{\rm (b)}" If $r(\pmb a)=\sigma_{k+1}$, then $f(N(\pmb a))=0$ and $r(N(\pmb a))=\sigma_{k+2}$.\endroster
\endproclaim 

\demo{Proof} It is clear that
$$\sum_{u=1}^{k+1}N(\pmb a)_u\sigma_u= \sum_{u=1}^{k}a_u\sigma_u+f(\pmb a)\sigma_{k+1}<n.$$We conclude that $N(\pmb a)$ is an eligible $(k+1)$-tuple. 
Suppose that $\pmb b$ is an eligible $j$-tuple with $\pmb a\gneq \pmb b\geq N(\pmb a)$.
Since $\pmb a \ge \pmb b\ge N(\pmb a)$ we have $j\ge k$ and $b_i=a_i$ for $i\le k$. As $\pmb a\gneq \pmb b$ we also have $j>k$. Now the inequality $\pmb b\ge N(\pmb a)$ implies $j=k+1$ and $b_{k+1}\ge f(\pmb a)$. Finally, the definition of $f(a)$ ensures $b_{k+1}\le f(\pmb a)$. Thus, $\pmb b=N(\pmb a)$. 
Assertion   (a) is established. 

The hypothesis of (b) yields 
$$\sigma_{k+1}=\sum_{u=1}^{k} a_u\sigma_u+(f(\pmb a)+1)\sigma_{k+1}-n+1;$$ hence,
$$n-1=\sum_{u=1}^{k} a_u\sigma_u+(f(\pmb a))\sigma_{k+1}=\sum_{u=1}^{k+1} N(\pmb a)_u\sigma_u.$$We now see that 
$f(N(\pmb a))=0$ and 
$$r(N(\pmb a))=\sum_{u=1}^{k+1} N(\pmb a)_u\sigma_u+(f(N(\pmb a))+1)\sigma_{k+2}-n+1=\sigma_{k+2}. \qed$$
\enddemo

\demo{Proof of Proposition \tref{:.O37.1*}} Once items (a) and (b) are shown, then Observation \tref{:.O28'} implies that the modules $\Cal E'_{\pmb a}/\Cal D'_{\pmb a}$ are Cohen-Macaulay, hence perfect. 

\medskip \flushpar (a) It suffices to show that $\Cal D'_{\pmb a}=\Cal D_{\pmb a}$. Let $\pmb b>\pmb a$ be the eligible tuple which is adjacent to $\pmb a$. It suffices to show that $\pmb b$ is eligible$'$. 
Suppose that $\pmb b$ is a $j$-tuple. If $\pmb b$ is not eligible$'$, then   $j<\ell -1$ and $r(\pmb b)=\sigma_{j+1}$. Now,
Lemma \tref{:sd42} shows   $r(\pmb a)=\sigma_{k+1}$ and this contradicts the hypothesis. 

\medskip \flushpar (b) Notice  that $k$ is necessarily equal to $\ell-1$. 
Identify the largest  
  non-negative integer $s$ for which there exists  an eligible   $(\ell-1-s)$-tuple $\pmb b$ with 
$\pmb a=N^{s}(\pmb b)$ and 
$r(\pmb b)=\sigma_{\ell-s}$. Let $j=\ell-1-s$. We know, from Lemma \tref{:sd42}, that
$$\pmb b>N(\pmb b)>N^2(\pmb b)>\dots>N^{s}(\pmb b)=\pmb a$$ are adjacent eligible neighbors and that if $0\le h\le s-1$, then $N^h(\pmb b)$ is not eligible$'$. Furthermore,  for each integer $h$, with $1\le h\le s$, we have $$f(N^h(\pmb b))=0\and r(N^h(\pmb b))=\sigma_{j+h+1}\tag\tnum{calc}$$ 
The module $\Cal E_{\pmb a}'/\Cal D_{\pmb b}$ is defined to be
$$\sum_{h=0}^sT^{N^h(\pmb b)} T_{j+1+h,1}^{f(N^h(\pmb b))}(T_{j+1+h,1},\dots, T_{j+1+h,r(N^h(\pmb b))})(S/\Cal D_{\pmb b}).$$ 
The calculations of (\tref{calc}) show that 
$$\Cal E_{\pmb a}'/\Cal D_{\pmb b}=T^{\pmb b}T_{j+1,1}^{f(\pmb b)} J(S/\Cal D_{\pmb b}),$$ where $J$ is generated by the entries in the first row of $\psi_{>j}$. We also know that $$T^{\pmb b}T_{j+1,1}^{f(\pmb b)}=T^{\pmb a}T_{k+1,1}^{f(\pmb a)}.$$ Furthermore,
 $\pmb\e$, which is defined to be $\e_{\pmb a}+f(\pmb a)\e_{k+1}$, is also equal to $\e_{\pmb b}+f(\pmb b)\e_{j+1}$.
Lemma \tref {:+.L37.!}(b) shows that $P_j$ annihilates $T^{\pmb b}T_{j+1,1}^{f(\pmb b)}J(S/\Cal D_{\pmb b})$. Apply  Lemma \tref{:L44.70} to the ideal $\Cal J=T_{j+1,1}^{f(\pmb b)}J$ to see that 
$$\split \Cal E_{\pmb a}'/\Cal D_{\pmb b}&{}=
T^{\pmb b}T_{j+1,1}^{f(\pmb b)}J(S/\Cal D_{\pmb b})\iso T_{j+1,1}^{f(\pmb b)}J(A/P_jA)(-\pmb\sigma\cdot\pmb \e_{\pmb b},0;-\pmb \e_{\pmb b})\\&{}\iso J(A/P_jA)(-\pmb\sigma\cdot\pmb \e,0;-\pmb \e).\endsplit$$
The final isomorphism holds because $T_{j+1,1}$ is a non-zero element in the domain $A/P_jA$.

Let $\pmb c$ be the eligible $i$-tuple for which $\pmb c>\pmb b$  are adjacent eligible neighbors. If $\pmb c$ is not eligible$'$, then   $i<\ell-1$ and $r(\pmb c)=\sigma_{i+1}$. Lemma \tref{:sd42}(a) then says that $\pmb b=N(\pmb c)$ and this contradicts the choice of $s$. Thus, $\pmb c$ is  eligible$'$, $\Cal D'_{\pmb a}=\Cal D_{\pmb b}$, and the proof is complete.  

\medskip \flushpar (c) Notice that $0^{\ell-1}$ is an eligible$'$-tuple. Let $\pmb b$ be the eligible $j$-tuple with $\pmb b>0^{\ell-1}$ and $\pmb b$ adjacent to $0^{\ell-1}$. It suffices to show that $\pmb b$ is eligible$'$. If $\pmb b$ is not eligible$'$, then $j<\ell-1$ and $r(\pmb b)=\sigma_{j+1}$. Lemma \tref{:sd42} then shows that $N(\pmb b)=0^{\ell-1}$; hence,
$$f(0^{\ell-1})=f(N(\pmb b))=0\and r(0^{\ell-1})=r(N(\pmb b))=\sigma_{\ell}.$$
The definition of $r$ now gives $\sigma_{\ell}=r(0^{\ell-1})=\sigma_{\ell}-n+1$; or $n=1$, which is a violation of the ambient hypotheses of Data \tref{:.data37}. 
 \qed
\enddemo

%%%%%%%%%%%%%%%% 
%% Section 6  %%
%%%%%%%%%%%%%%%%
%\newpage
\bigpagebreak
\SectionNumber=6\tNumber=1
\flushpar{\bf \number\SectionNumber.\quad  Symbolic Rees Algebra.}

\medskip

Retain the notation of (\tref{:.data37}).
\proclaim{Proposition \tnum{P6.1}} The symbolic Rees algebra
$$\Cal R_{s}(K)=\bigoplus_{n\ge 0} K^{(n)}$$is finitely generated as an $A$-algebra.
\endproclaim

\demo{Proof} View $\Cal R_{s}(K)$ as the subring of the polynomial ring $A[u]$ which is generated by 
$$\bigcup\limits_{n=1}^{\infty} \{\theta u^n\mid \theta\in K^{(n)}\}.$$ Let $\Cal S$ be the following set of elements of $\Cal R_{s}(K)$ 
$$\Cal S=\{T_{i,j}u^k\mid 1\le i\le \ell,\quad 1\le j\le \sigma_i,\and 1\le k\le \sigma_i+1-j\}.$$
We prove that $\Cal R_{s}(K)$ is generated as an $A$-algebra by $\Cal S$. Suppose that $\theta$ is a generator of $K^{(n)}$. Then there is an eligible $k$-tuple $\pmb a$ with $\theta= T^{\pmb a}T_{k+1,1}^{f(\pmb a)}T_{k+1,j}$. We have 
$$\sum\limits_{u=1}^k a_u\sigma_u+f(\pmb a)\sigma_{k+1}<n\and 1\le j\le r(\pmb a)\le \sigma_{k+1}.$$ Thus,
$$\theta u^n=\prod\limits_{i=1}^k\(T_{i,1}u^{\sigma_i}\)^{a_i}\( T_{k+1,1}u^{\sigma_{k+1}}\)^{f(\pmb a)}T_{k+1,j}u^{\sigma_{k+1}+1-r(\pmb a)}\in A[\Cal S].
\qed$$
\enddemo
 
%%%%%%%%%%%%%%%% 
%% References %%
%%%%%%%%%%%%%%%%
%\newpage
\medskip

\flushpar {\bf Acknowledgment.}
This work was conducted while the first author was on sabbatical at Purdue University and later was a Visiting Professor at the University of Notre Dame. He appreciates the sabbatical from the University of South Carolina and the hospitality he received at Purdue University and the University of Notre Dame.

\enddocument